\documentclass[12pt]{amsart}
\usepackage{latexsym, amssymb, amsmath, amscd}
 \usepackage{eucal}

    \newtheorem{rema}{Remark}[section]
    \newtheorem{propo}[rema]{Proposition}

   \newtheorem{theo}[rema]{Theorem}
   \newtheorem{def-theo}[rema]{Definition-Theorem}
 
   \newtheorem{defi}[rema]{Definition}
    \newtheorem{lemma}[rema]{Lemma}
    \newtheorem{corol}[rema]{Corollary}
     
  \newtheorem{rmk}[rema]{Remark}

	\newcommand{\nno}{\nonumber}

	\newcommand{\p}{\partial}

 \newcommand{\pf}{{\it Proof:}\hspace{2ex}}
 
 \newcommand{\epfv}{\hspace{1em}$\Box$\vspace{1em}}

\newcommand{\bZ}{{\mathbb Z}}

\newcommand{\bQ}{{\mathbb Q}}
\newcommand{\bN}{{\mathbb N}}
\newcommand{\bT}{{\mathbb T}}
\newcommand{\bbT}{{\bar{\mathbb T}}}

\newcommand{\cT}{{\mathcal T}}
\newcommand{\cTW}{\cT_W}
\newcommand{\cP}{{\mathcal P}}
\newcommand{\cA}{{\mathcal A}}
\newcommand{\cS}{{\mathcal S}}

\newcommand{\cD}{{\mathcal D}}
\newcommand{\cH}{{\mathcal H}}

\newcommand{\cNcs}{{${\mathcal N}$CS} }
\newcommand{\cNsf}{{{\mathcal N}Sym}}
\newcommand{\cQf}{{{\mathcal Q}Sym}}
\newcommand{\Oft}{ \Omega_{F_t}  }
\newcommand{\OTW}{{\Omega_\bT^W} }
\newcommand{\cSft}{ {\mathcal S}_{F_t} }
\newcommand{\cAft}{ {\mathcal A}_{F_t} }

\newcommand{\cDz}{{\mathcal D \langle  z \rangle}}
\newcommand{\cDaz}{{{\mathcal D}^{[\alpha]} \langle  z \rangle}}

\newcommand{\cDzz}{{\mathcal D \langle \langle z \rangle\rangle}}
\newcommand{\cDkzz}{{\mathcal D_K \langle \langle z \rangle\rangle}}
\newcommand{\cDtzz}{{\mathcal D_t \langle \langle z \rangle\rangle}}
\newcommand{\cDrzz}{{\mathcal D er\langle \langle z \rangle\rangle}}
\newcommand{\cDrtzz}{{\mathcal D er_t \langle \langle z \rangle\rangle}}
\newcommand{\cDrkzz}{{\mathcal D er_K \langle \langle z \rangle\rangle}}

\newcommand{\cDrazz}{{\cD er^{[\alpha]}\langle \langle z \rangle \rangle}} 
 
\newcommand{\cDazz}{{\cD^{[\alpha]}\langle \langle z \rangle \rangle}}

\newcommand{\cDrtazz}{{\cD er^{[\alpha]}_t\langle \langle z \rangle \rangle}} 
 
\newcommand{\cDtazz}{{\cD^{[\alpha]}_t\langle \langle z \rangle \rangle}}

\newcommand{\ataz}{{\mathbb A_t^{[\alpha]}\langle \langle z\rangle\rangle}}
\newcommand{\gtaz}{{\mathbb G_t^{[\alpha]}\langle \langle z\rangle\rangle}}

\newcommand{\kz}{{K\langle z \rangle}}
\newcommand{\kzz}{{K\langle \langle z \rangle\rangle}}
\newcommand{\ktz}{{K[t]\langle z \rangle}}
\newcommand{\ktzz}{{K[t]\langle \langle z \rangle\rangle}}
\newcommand{\kttzz}{{K[[t]]\langle \langle z \rangle\rangle}}

\newcommand{\BQ}{\begin{eqnarray}}
\newcommand{\EQ}{\end{eqnarray}}
\newcommand{\BQn}{\begin{eqnarray*}}
\newcommand{\EQn}{\end{eqnarray*}}

\newcommand{\lb}{\left[}
\newcommand{\rb}{\right]}
\newcommand{\lp}{\left(}
\newcommand{\rp}{\right)}

\newcommand{\fr}{\frac}
\newcommand{\pz}{\frac{\p}{\p z}}
\newcommand{\wtilde}{\widetilde}

\title[Noncommutative Symmetric   Systems]
{\cNcs  Systems over Differential Operator Algebras 
and the Grossman-Larson Hopf Algebras
 of Labeled Rooted Trees}

    \author{Wenhua Zhao}      
    \begin{document}

\begin{abstract}
Let $K$ be any unital commutative $\bQ$-algebra and $W$ any non-empty 
subset of $\bN^+$. Let $z=(z_1, \ldots, z_n)$ 
be commutative or noncommutative free variables and $t$ a 
formal central parameter. 
Let $\cDazz$ $(\alpha\geq 1)$ be the unital algebra 
generated by the differential operators of $\kzz$ 
which increase the degree in $z$ by at least $\alpha-1$ 
and $ \ataz $ the group of automorphisms $F_t(z)=z-H_t(z)$ 
of $\kttzz$ with $o(H_t(z))\geq \alpha$ and $H_{t=0}(z)=0$.
First, we study a connection of the \cNcs  systems 
$\Omega_{F_t}$ $(F_t\in \ataz)$
(\cite{GTS-I}, \cite{GTS-II}) 
over the differential operators algebra 
$\cDazz$ and the \cNcs  system $\Omega_\bT^W$ (\cite{GTS-IV}) 
over the Grossman-Larson 
Hopf algebra $\cH_{GL}^W$ (\cite{GL}, \cite{F1}, \cite{F2})
of $W$-labeled rooted trees. 
We construct a Hopf algebra homomorphism 
$\mathcal A_{F_t}: \cH_{GL}^W \to \cDazz$ 
$(F_t\in \ataz)$ such that 
$\mathcal A_{F_t}^{\times 5}(\Omega_\bT^W)
=\Omega_{F_t}$. Secondly, 
we generalize the tree expansion formulas for 
the inverse map (\cite{BCW}, 
\cite{Wr3}), the D-Log and the formal flow (\cite{WZ}) 
of $F_t$ in the commutative case 
to the noncommutative case.
Thirdly, we prove the injectivity of  
the specialization $\cT:{\mathcal N}Sym \to \cH_{GL}^{\bN^+}$ 
(\cite{GTS-IV}) of NCSF's 
(noncommutative symmetric functions) (\cite{G-T}). 
Finally, we show  
the family of the specializations
$\cS_{F_t}$ of NCSF's with all $n\geq 1$ 
and the polynomial automorphisms $F_t=z-H_t(z)$
with $H_t(z)$ homogeneous and 
the Jacobian matrix $JH_t$ strictly lower triangular
can distinguish any two different NCSF's. 
The graded dualized versions of the main results 
above are also discussed. 
\end{abstract}

   

 \bibliographystyle{alpha}
    \maketitle


\renewcommand{\theequation}{\thesection.\arabic{equation}}
\renewcommand{\therema}{\thesection.\arabic{rema}}
\setcounter{equation}{0}
\setcounter{rema}{0}
\setcounter{section}{0}

\section{\bf Introduction}
 \label{S1}

Let $K$ be any unital commutative  
$\bQ$-algebra and $A$ a 
unital associative 
but not necessarily 
commutative $K$-algebra.
A {\it \cNcs $($non-commutative symmetric$)$  system} over $A$ 
(see Definition \ref{Main-Def}) 
by definition is a $5$-tuple 
$\Omega\in A[[t]]^{\times 5}$ 
which satisfies the defining equations 
(see Eqs.\,$(\ref{UE-0})$--$(\ref{UE-4})$) 
of the NCSF's (noncommutative symmetric functions) 
first introduced and studied 
in the seminal paper \cite{G-T}. 
When the base algebra 
$K$ is clear in the context,
the ordered pair $(A, \Omega)$ 
is also called a {\it \cNcs system}.
In some sense, a \cNcs  system over 
an associative $K$-algebra can be viewed 
as a system of analogs in $A$
of the NCSF's defined by 
Eqs.\,$(\ref{UE-0})$--$(\ref{UE-4})$.
For some general discussions on 
the \cNcs  systems, see \cite{GTS-I}.
For more studies on NCSF's, 
see \cite{T}, \cite{NCSF-II}, 
\cite{NCSF-III}, \cite{NCSF-IV}, 
\cite{NCSF-V} and \cite{NCSF-VI}.
One immediate but probably the most 
important example of the \cNcs systems 
is $(\cNsf, \Pi)$ formed 
by the generating functions of 
the NCSF's defined in \cite{G-T}
by Eqs.\,$(\ref{UE-0})$--$(\ref{UE-4})$ 
over the free $K$-algebra $\cNsf$ of NCSF's  
(see Section \ref{S2.1}). 
It serves as the universal \cNcs  system 
over all associative $K$-algebra 
(see Theorem \ref{Universal}).  
More precisely, for any \cNcs  system $(A, \Omega)$, 
there exists a unique $K$-algebra homomorphism 
$\cS: \cNsf \to A$ such that 
$\cS^{\times 5}(\Pi) = \Omega$ 
(here we have extended the homomorphism
$\cS$ to $\cS: \cNsf[[t]] \to A[[t]]$ 
by the base extension). 
In \cite{GTS-IV} and \cite{GTS-II}, 
some families of \cNcs  systems 
over differential operator algebras and 
the Grossman-Larson Hopf algebra 
(\cite{GL}, \cite{F1}, \cite{F2}) 
of labeled rooted trees have 
been constructed, respectively. 
Consequently, by the universal property 
of the \cNcs  system $(\cNsf, \Pi)$, 
one obtains two families 
of specializations of NCSF's 
by differential operators 
and labeled rooted trees 
(see Sections \ref{S2.2} 
and \ref{S2.3} for a brief 
review of the results above).

In the first part of this paper, we study a
connection of the \cNcs  systems 
in \cite{GTS-IV} over 
the Grossman-Larson Hopf algebras 
of labeled rooted trees and 
the \cNcs  systems in \cite{GTS-II} over 
differential operator algebras. 
We construct a Hopf algebra homomorphism 
from the former algebra to the later 
which maps the former \cNcs  system 
to the later. 
In the second part of this paper,  
we first apply the connection above to 
derive tree expansion formulas for 
the D-Log's, the formal flows and 
the inverse maps 
of formal automorphisms 
in commutative or noncommutative 
variables, which generalize the 
tree expansion formulas obtained in \cite{BCW}, 
\cite{Wr3} and \cite{WZ} for commutative 
variables to the noncommutative case. 
Secondly, by combining the connection above 
with some results obtained in \cite{GTS-II} 
and \cite{GTS-IV}, we prove more properties 
for the specializations of NCSF's 
by differential operator 
and labeled rooted trees 
obtained in \cite{GTS-II} 
and \cite{GTS-IV}, respectively.

To be more precise, let $z=(z_1, \ldots, z_n)$ 
be commutative or noncommutative free variables 
and $t$ a formal central parameter. 
Denote uniformly for both commutative 
and noncommutative variables $z$
by $\kzz$ (resp.\,\,$\kz$) 
the formal power series (resp.\,\,polynomial) 
algebra of $z$ over $K$.
For any $\alpha \geq 1$, 
let $\cDazz$ (resp.\,\,$\cDaz$) 
be the unital algebra generated by 
the differential operators of $\kzz$ (resp.\,\,$\kz$)
which increase the degree in $z$ by at least 
$\alpha-1$ and $ \ataz $ the group of automorphisms $F_t(z)=z-H_t(z)$ 
of $\kttzz$ with $o(H_t(z))\geq \alpha$ and $H_{t=0}(z)=0$.
In \cite{GTS-II}, associated with each automorphism 
$F_t\in \ataz$, a \cNcs  system (\cite{GTS-I}) $\Oft$ over the 
differential operator algebra $\cDazz$ 
has been constructed. Consequently, 
by the universal property 
of the \cNcs  system $(\cNsf, \Pi)$, 
one obtains a families 
of specializations $\cSft: \cNsf \to\cDazz$ 
of NCSF's by differential operators. 
In \cite{GTS-IV}, for any non-empty 
$W\subseteq \bN^+$, a \cNcs  system $\OTW$ over 
the Grossman-Larson Hopf algebra 
$\cH_{GL}^W$ were given explicitly.
Hence, one also gets a specialization
$\cTW: \cNsf\to \cH_{GL}^W$
of NCSF's by $W$-labeled rooted trees. 

In the first part of this paper, 
for any fixed $\alpha \geq 1$, $\emptyset\neq W\subseteq \bN^+$
and $F_t(z)=z-H_t(z)\in \ataz$  such that 
$H_t(z)$ can be written as 
$\sum_{m\in W} t^m H_{[m]}(z)$ for some 
$H_{[m]}(z)\in \kzz^{\times n}$ $(m\in W)$,
we construct a Hopf algebra
homomorphism $\cAft: \cH_{GL}^W \to \cDazz$
such that $\cAft^{\times 5}(\OTW)=\Oft$ 
(see Theorem \ref{S3-Main-1} and  \ref{S3-Main-2}).
Furthermore, we also show in 
Proposition \ref{CD-2} 
that the specializations 
$\cSft: \cNsf \to \cDazz$ and
$\cTW: \cNsf \to \cH_{GL}^W$ of NCSF's 
are connected by $\cSft=\cTW \circ\cAft$.
Note that, it has been shown in \cite{GTS-II} 
that $\cSft$ is a graded Hopf algebra 
homomorphism from $\cNsf$ to 
the Hopf subalgebra 
$\cDaz\subset \cDazz$
iff $F_t\in \ataz$
has the form $F_t(z)=t^{-1} F(tz)$ 
for some automorphism 
$F(z)$ of $\kzz$. 
By taking the graded duals 
of the results above, 
we get the corresponding commutative diagrams
(see Proposition \ref{CD-3}) 
for the Hopf algebras $\cQf$ 
of quasi-symmetric functions 
(\cite{Ge}, \cite{MR}, \cite{St2}), 
the Connes-Kreimer Hopf algebra $\cH_{CK}^W$ 
of $W$-labeled rooted forests 
(\cite{CM}, \cite{Kr}, \cite{CK}, \cite{F1}, \cite{F2}), 
and the graded dual $\cD^{[\alpha]}\langle z\rangle^*$ of 
differential operator algebra $\cD^{[\alpha]}\langle z\rangle$ 
of $\kz$. 

In the second part of this paper,  
by applying the Hopf algebra 
homomorphism $\cAft: \cH_{GL}^W \to \cDazz$ 
described above, we first derive tree expansion formulas 
for the D-Log, the formal flow and the inverse map of $F_t$ 
(see Theorem \ref{S4-Main} and Corollaries \ref{M-powers}, 
\ref{S4-Main-C}).  Since the proofs given here 
do not depend on the commutativity 
of the free variables $z$, 
the formulas derived here can be viewed 
as some natural generalizations 
to the noncommutative case of 
the tree expansion formulas derived 
in \cite{BCW}, \cite{Wr3} and \cite{WZ}  
in the commutative case.
Finally, we apply the Hopf algebra 
homomorphism $\cAft: \cH_{GL}^W \to \cDazz$ 
above combining with some results already 
obtained in \cite{GTS-II} 
and \cite{GTS-IV} to study more properties 
of the specializations $\cSft: \cNsf \to \cDazz$ 
and $\cT_W: \cNsf \to \cH_{GL}^W$ of NCSF's. 
In Theorem \ref{NCSF-Injc-Trees},
we show that, when $W=\bN^+$, 
the specialization
$\cTW: \cNsf \to \cH_{GL}^W$ actually 
embeds $\cNsf$ into $\cH_{GL}^W$ 
as a graded $K$-Hopf subalgebra. 
By taking the graded duals, we get a surjective 
graded Hopf algebra homomorphism 
$\cT_W^*: \cH_{CK}^W \to \cQf$
from the Connes-Kreimer Hopf algebra $\cH_{CK}^W$ 
(\cite{CM}, \cite{Kr}, \cite{CK}, \cite{F1}, \cite{F2})
of $W$-labeled rooted forests 
onto the Hopf algebra $\cQf$  
of quasi-symmetric functions 
(see Proposition \ref{Onto}).
In Theorem \ref{StabInjc-best}, we show  
the family of the differential operator 
specializations $\cS_{F_t}$ of NCSF's 
with all $n\geq 1$  and 
polynomial automorphisms $F_t=z-H_t(z) \in \ataz$
such that, $H_t(z)$ is homogeneous and 
the Jacobian matrix $JH$ 
is strictly lower triangular,
can distinguish any 
two different NCSF's. 

Considering the length of paper, 
we give a more detailed arrangement 
of the paper as follows.

In Section \ref{S2}, we mainly fix some
necessary notations 
and recall some main results of \cite{GTS-II} 
and \cite{GTS-IV} that will be needed 
throughout this paper. In Subsection \ref{S2.1}, 
we recall the definition of 
general \cNcs  systems and the universal \cNcs  system 
$(\cNsf, \Pi)$ from NCSF's. 
In Subsection \ref{S2.2}, we recall the \cNcs  systems 
$\Oft$ constructed in \cite{GTS-II} over 
the differential operator algebras $\cDazz$
and the resulted differential 
operator specializations 
$\cSft: \cNsf \to \cDazz$ 
of NCSF's. 
In Subsection \ref{S2.3}, 
we recall the definition of 
the \cNcs  systems $\cH_{GL}^W$ 
$(\emptyset \neq W \subseteq \bN^+)$ 
constructed in \cite{GTS-IV} over the Grossman-Larson 
Hopf algebra of $W$-labeled rooted trees.
In Section \ref{S3}, for any 
$F_t(z)=z-H_t(z)\in \ataz$ and any
non-empty $W\subseteq \bN^+$ such that 
$H_t(z)=\sum_{m\in W} t^m H_{[m]}(z)$ for some 
$H_{[m]}(z)\in \kzz^{\times n}$ $(m\in W)$,
we constructed a Hopf algebra homomorphisms 
$\cAft:\cH_{GL}^W \to \cDazz$ 
such that $\cAft^{\times 5}(\OTW)=\Oft$ 
(see Theorem \ref{S3-Main-1} and \ref{S3-Main-2}).
Furthermore, we also show in 
Proposition \ref{CD-2} that the specializations 
$\cSft: \cNsf \to \cDazz$ and
 $\cTW: \cNsf \to \cH_{GL}^W$ of NCSF's 
are connected by $\cSft=\cTW \circ\cAft$.
The graded dualized versions of the main results
above are also discussed.
In Section \ref{S4}, by applying the Hopf algebra homomorphism 
$\cAft: \cH_{GL}^W \to \cDazz$ constructed 
in the previous section, we derive tree expansion formulas 
for the D-Log, formal flow and the inverse map of $F_t$ 
(see Theorems \ref{S4-Main} and Corollaries \ref{M-powers}, 
\ref{S4-Main-C}).  
In Section \ref{S5}, 
we study more properties 
of the specializations $\cSft: \cNsf \to \cDazz$ 
and $\cTW: \cNsf \to \cH_{GL}^W$ of NCSF's. First, 
we show in Theorem \ref{NCSF-Injc-Trees} that, 
when $W=\bN^+$, the specialization
$\cTW: \cNsf \to \cH_{GL}^W$ is actually 
an injective graded $K$-Hopf 
algebra homomorphism. By taking the graded duals, 
we get a surjective graded Hopf algebra
homomorphism $\cTW^*: \cH_{CK}^W\to \cQf$
(see Proposition \ref{Onto}).  Secondly,
we show in Theorem \ref{StabInjc-best} that
the family of the specializations
$\cS_{F_t}$ of NCSF's with all $n\geq 1$ 
and the polynomial automorphisms $F_t=z-H_t(z)$
with $H_t(z)$ homogeneous and 
the Jacobian matrix $JH_t$ 
strictly lower triangular
can distinguish 
any two different NCSF's.

\renewcommand{\theequation}{\thesection.\arabic{equation}}
\renewcommand{\therema}{\thesection.\arabic{rema}}
\setcounter{equation}{0}
\setcounter{rema}{0}

\section{\bf \cNcs  Systems} \label{S2}

Let $K$ be any unital commutative $\bQ$-algebra and 
$A$ any unital associative but not necessarily commutative 
$K$-algebra. Let $t$ be a formal central parameter, 
i.e. it commutes with all elements of $A$, and $A[[t]]$ 
the $K$-algebra of formal power series 
in $t$ with coefficients in $A$. 
First let us recall the following notion formulated 
in \cite{GTS-I}.

\begin{defi} \label{Main-Def}
For any  unital associative $K$-algebra $A$, a $5$-tuple $\Omega=$ 
$( f(t)$, $g(t)$, $d\,(t)$, $h(t)$, $m(t) ) 
\in A[[t]]^{\times 5}$ is said 
to be a {\it \cNcs $($noncommutative symmetric$)$  system}
over $A$ if the following equations are satisfied.
\allowdisplaybreaks{
\begin{align}
&f(0)=1 \label{UE-0}\\
& f(-t)  g(t)=g(t)f (-t)=1, \label{UE-1}   \\
& e^{d\,(t)} = g(t), \label{UE-2} \\
& \frac {d g(t)} {d t}= g(t) h(t), \label{UE-3}\\ 
& \frac {d g(t)}{d t} =  m(t) g(t).\label{UE-4}
\end{align}}
\end{defi} 

When the base algebra $K$ is clear in the context, we also call 
the ordered pair $(A, \Omega)$ a {\it \cNcs  system}. 
Since \cNcs  systems often come from generating functions 
of certain elements of $A$ that are under concern, 
the components of $\Omega$ will also be refereed as 
the {\it generating functions} of their coefficients. 

Since all $K$-algebras that we are going to work on in this paper
$K$-Hopf algebras (\cite{Abe}, \cite{Knu}, \cite{Mon}), 
the following result proved in \cite{GTS-I} 
later will be useful to our later arguments.

\begin{propo}\label{bialg-case}
Let $(A, \Omega)$ be a \cNcs  system as above. 
Suppose $A$ is further  a $K$-bialgebra. 
T\begin{enumerate}
\item[$(a)$] The coefficients of $f(t)$ form a divided power series of $A$.

\item[$(b)$] The coefficients of $g(t)$ form a divided power series of $A$.

\item[$(c)$] One $(\text{hence also all})$ 
of \,$d(t)$, $h(t)$ and $m(t)$ has all 
its coefficients primitive in $A$.
\end{enumerate}
\end{propo}

In this section, we briefly recall in 
Subsection \ref{S2.1} 
the \cNcs  system $(\cNsf, \Pi)$ formed by 
generating functions of some of the NCSF's defined in 
\cite{G-T} 
and its universal property 
(see Theorem \ref{Universal}).  
In Subsection \ref{S2.2}, we recall
the \cNcs  systems constructed in \cite{GTS-II}
over differential operator algebras.
Finally, in Subsection \ref{S2.3}, we recall
the \cNcs  systems constructed in \cite{GTS-IV}
over the Grossman-Larson 
Hopf algebra of labeled rooted trees.

\subsection{The Universal \cNcs  System from
Noncommutative Symmetric Functions}\label{S2.1}

Let $\Lambda=\{ \Lambda_m\,|\, m\geq 1\}$ 
be a sequence of noncommutative 
free variables and $\cNsf$ 
the free associative algebra 
generated by 
$\Lambda$ over $K$.  For convenience, 
we also set $\Lambda_0=1$.
We denote by
$\lambda (t)$ the generating function of 
$\Lambda_m$ $(m\geq 0)$, i.e. we set
\begin{align}
\lambda (t):= \sum_{m\geq 0} t^m \Lambda_m 
=1+\sum_{k\geq 1} t^m \Lambda_m.
\end{align}

In the theory of NCSF's (\cite{G-T}), 
$\Lambda_m$ $(m\geq 0)$ is 
the noncommutative analog 
of the $m^{th}$ classical (commutative) 
elementary symmetric function 
and is called the {\it $m^{th}$ 
$(\text{noncommutative})$ 
elementary symmetric function.}

To define some other NCSF's, we consider 
Eqs.\,$(\ref{UE-1})$--$(\ref{UE-4})$ 
over the free $K$-algebra $\cNsf$
with $f(t)=\lambda(t)$. The 
solutions for $g(t)$, $d\,(t)$, 
$h(t)$, $m(t)$ exist and are unique, 
whose coefficients will be the NCSF's 
that we are going to define.
Following the notation in \cite{G-T} 
and \cite{G-T}, we denote the resulted 
$5$-tuple by 
\begin{align}
\Pi:= (\lambda(t),\, \sigma(t),\, \Phi(t),\, \psi(t),\, \xi(t))
\end{align}
and write the last 
four generating functions of 
$\Pi$ explicitly as follows.

\allowdisplaybreaks{
\begin{align}
\sigma (t)&=\sum_{m\geq 0} t^m S_m,  \label{lambda(t)} \\
\Phi (t)&=\sum_{m\geq 1} t^m \frac{\Phi_m}m  \label{Phi(t)}\\
\psi (t)&=\sum_{m\geq 1} t^{m-1} \Psi_m, \label{psi(t)}\\
\xi (t)&=\sum_{m\geq 1} t^{m-1} \Xi_m.\label{xi(t)}
\end{align}}

Following \cite{G-T},
we call $S_m$ ($m\geq 1$) the 
{\it $m^{th}$ $(\text{noncommutative})$ complete 
homogeneous symmetric function} and
$\Phi_m $ (resp.\,\,$\Psi_m$) 
the {\it $m^{th}$ power sum symmetric function 
of the second $($resp.\,\,first$)$ kind}. 
Following \cite{GTS-I}, 
we call $\Xi_m \in \cNsf$ $(m\geq 1)$ 
the {\it $m^{th}$ $(\text{noncommutative})$ 
power sum symmetric function of the third kind}.

Next, let us recall the following graded 
$K$-Hopf algebra structure 
of $\cNsf$. It has been shown in 
\cite{G-T} that $\cNsf$ is the universal enveloping algebra 
of the free Lie algebra generated 
by $\Psi_m$ $(m\geq 1)$. Hence, it has a Hopf  
$K$-algebra structure as all other universal enveloping algebras 
of Lie algebras do. Its co-unit $\epsilon:\cNsf \to K$,
 co-product $\Delta$ and 
 antipode $S$ are uniquely determined by 
\begin{align}
\epsilon (\Psi_m)&=0, \label{counit} \\
\Delta (\Psi_m) &=1\otimes \Psi_m +\Psi_m\otimes 1, \label{coprod}\\
S(\Psi_m) & =-\Psi_m,\label{antipode}
\end{align}
for any $m\geq 1$. 

Next, we introduce the {\it weight} of NCSF's 
by setting the weight of 
any monomial $\Lambda_{m_1}^{i_1} 
\Lambda_{m_2}^{i_2} \cdots \Phi_{m_k}^{i_k}$
to be $\sum_{j=1}^k i_j m_j$. 
For any $m\geq 0$, we denote by $\cNsf_{[m]}$ 
the vector subspace of $\cNsf$ spanned 
by the monomials of $\Lambda$ 
of weight $m$. Then it is easy to see that 
\begin{align}\label{Grading-cNsf}
\cNsf=\bigoplus_{m\geq 0} \cNsf_{[m]}, 
\end{align}
which provides a grading for $\cNsf$. 

Note that, it has been shown in \cite{G-T}, 
for any $m\geq 1$, the NCSF's 
$S_m, \Phi_m, \Psi_m \in  \cNsf_{[m]}$. 
By the facts above and 
Eqs.\,(\ref{counit})--(\ref{antipode}), 
it is also easy to check that, 
with the grading given in Eq.\,(\ref{Grading-cNsf}), 
$\cNsf$ forms a graded $K$-Hopf algebra. 
Its graded dual is given 
by the space $\cQf$ of quasi-symmetric functions, 
which were first introduced by I. Gessel \cite{Ge} 
(also see \cite{MR} and \cite{St2} for more discussions).

Now we come back to our discussions on the \cNcs  systems. 
From the definitions of the NCSF's above, 
we see that $(\cNsf, \Pi)$ obviously forms a \cNcs  system.
More importantly, as shown in Theorem $2.1$ in \cite{GTS-I}, 
we have the following important theorem on 
the \cNcs  system $(\cNsf, \Pi)$. 

\begin{theo}\label{Universal}
Let $A$ be a $K$-algebra and $\Omega$ 
a \cNcs  system over $A$. Then, 

$(a)$ There exists a unique $K$-algebra homomorphism 
$\cS: \cNsf\to A$ such that 
$\cS^{\times 5} (\Pi)=\Omega$.

$(b)$ If $A$ is further  a $K$-bialgebra $($resp.\,\,$K$-Hopf algebra$)$ 
and one of the equivalent statements in Proposition \ref{bialg-case} 
holds for the \cNcs  system $\Omega$, then $\cS: \cNsf\to A$ is also 
a homomorphism of $K$-bialgebras $($resp.\,\,$K$-Hopf algebras$)$.
\end{theo}

\begin{rmk}\label{Motivation}
By applying the similar arguments as in the 
proof of Theorem \ref{Universal}, 
or simply taking the quotient over 
the two-sided ideal generated by the commutators 
of $\Lambda_m$'s, it is easy to see that, 
over the category of commutative $K$-algebras, 
the universal \cNcs  system 
is given by the generating functions of
the corresponding classical 
$($commutative$)$ symmetric functions \cite{M}.
\end{rmk}

\subsection{\cNcs  Systems over Differential Operator Algebras}
\label{S2.2}

In this subsection, we briefly recall 
the \cNcs  systems constructed 
in \cite{GTS-II} over 
the differential operator 
algebras in commutative 
or noncommutative free variables. 
The construction of this \cNcs is 
mainly motivated by the studies 
in \cite{BurgersEq} and \cite{NC-IVP}
on the deformations of formal analytic 
maps and their applications 
to the inversion problem 
(\cite{BCW}, \cite{E}, 
\cite{GTS-III}).

First, let us fix the following notation.

Let $K$ be any unital 
commutative $\bQ$-algebra as before
and $z=(z_1, z_2, ... , z_n)$ commutative 
or noncommutative free 
variables.\footnote{Since most of the results 
as well as 
their proofs in this paper 
do not depend on the commutativity 
of the free variables $z$,  
we will not distinguish the commutative 
and the noncommutative case, 
unless stated otherwise, 
and adapt the notations 
for noncommutative variables 
uniformly for the both cases.}
Let $t$ be a formal central parameter, 
i.e. it commutes with $z$ and elements of $K$.
We denote by $\kzz$ and 
$\kttzz$ the $K$-algebras of 
formal power series in $z$ over  
$K$ and $K[[t]]$, respectively.
We denote by 
$\cDrkzz$ or $\cDrzz$, when the base algebra
$K$ is clear from the context, 
the set of all $K$-derivations 
of $\kzz$.
The unital subalgebra of 
$\text{End}_k(\kzz)$
generated by all
$K$-derivations 
of $\kzz$ will be denoted by 
$\cDkzz$ or $\cDzz$. 
Elements of $\cDkzz$ will be called 
{\it $(\text{formal})$ 
differential operators} in the
commutative and 
noncommutative variables $z$.

For any $\alpha\geq 1 $, 
we denote by $\cDrazz$ 
the set of the $K$-derivations 
of $\kzz$ which increase 
the degree in $z$ by 
at least $\alpha-1$. 
The unital subalgebra of 
$\cDzz$ generated by elements of 
$\cDrazz$ will be denoted by $\cDazz$. 
Note that, by the definitions above,
the operators of scalar multiplications
are also in $\cDzz$ and $\cDazz$.
When the base algebra is $K[[t]]$ 
instead of $K$ itself,
the notation 
$\cDrzz$, $\cDzz$, $\cDrazz$ 
and $\cDazz$ will be denoted by  
$\cDrtzz$, $\cDtzz$, $\cDrtazz$ 
and $\cDtazz$, respectively.
For example,  $\cDrtazz$ stands for 
the set of all $K[[t]]$-derivations of $\kttzz$ 
which increase the degree in $z$ 
by at least $\alpha-1$. 
Note that, $\cDrtazz=\cDrazz[[t]]$ and 
$\cDtazz=\cDazz[[t]]$.

For any $1\leq i\leq n$ and $u(z)\in \kzz$, 
we denote by $\lb u(z) \fr \p{\p z_i}\rb $
the $K$-derivation which maps $z_i$ to $u(z)$ and $z_j$ to $0$ 
for any $j\neq i$. 
For any $\Vec{u}=(u_1, u_2, \cdots, u_n)\in \kzz^{\times n}$, 
we set 
\begin{align}\label{Upz}
[\Vec{u}\pz]:=\sum_{i=1}^n [u_i \fr\p{\p z_i}].  
\end{align}

Note that, in the noncommutative case, 
we in general do {\bf not} have
$\lb u(z) \fr \p{\p z_i}\rb  g(z)  = u(z)  
\fr {\p g}{\p z_i}$ for all 
$u(z), g(z)\in \kzz$. This is the reason that 
we put a bracket $[\cdot]$ in the notation above 
for the $K$-derivations.
With the notation above, 
it is easy to see that any 
$K$-derivations $\delta$ 
of $\kzz$ can be written uniquely  as 
$\sum_{i=1}^n \lb f_i(z)\fr\p{\p z_i}\rb$ 
with $f_i(z)=\delta z_i\in \kzz$ 
$(1\leq i\leq n)$.

Note that, the differential operator algebra
$\cDazz$ $(\alpha\geq 1)$, as the universal 
enveloping algebra of
Lie algebra $\cDrazz$ with the commutator bracket, 
has a Hopf algebra structure as all 
other enveloping algebras of 
Lie algebras do.
In particular,
Its coproduct $\Delta$, antipode $S$ and co-unit $\epsilon$
are uniquely determined respectively by the properties
\begin{align}
\Delta(\delta)&= 1\otimes\delta+\delta\otimes 1,\label{Coprd-delta} \\
S(\delta)&=-\delta, \label{antipd-delta}\\
\epsilon (\delta) &=\delta \cdot 1, \label{Counit-delta}
\end{align}
for any $\delta \in \cDrzz$.

For any $\alpha\geq 1 $, let $\ataz$ 
be the set of all the automorphism 
$F_t(z)$ of $\kttzz$ over $K[[t]]$, 
which have the form
$F(z)=z-H_t(z)$ for some 
$H_t(z)\in \kttzz^{\times n}$ 
with $o(H_t(z))\geq \alpha$ and $H_{t=0}(z)=0$. 
Note that, for any $F_t \in \ataz$ as above, 
its inverse map $G_t:=F_t^{-1}$
can always be written uniquely 
as $G_t(z)=z+M_t(z)$ \label{Mt(z)}
for some $M_t(z)\in \kttzz^{\times n}$ 
with $o(M_t(z))\geq \alpha$ 
and $M_{t=0}(z)=0$. 

Now we recall the \cNcs  systems constructed 
in \cite{GTS-II} over 
the differential operator algebras 
$\cDazz$ $(\alpha\geq 1)$. 

We fix an $\alpha \geq 1$ and 
an arbitrary $F_t\in \ataz$. 
We will always 
let $H_t(z)$, $G_t(z)$ 
and $M_t(z)$ be determined 
as above. The \cNcs  system
\begin{align}\label{Def-Omega-Ft}
\Omega_{F_t}=(f(t), \, g(t),\, d(t),\, h(t),\, m(t)) \in \cDazz[[t]]^{\times 5}.
\end{align}
are determined as follows.

The last two components are given directly as

\begin{align}
h(t)&:=\left[\frac{\p M_t}{\p t}(F_t) \pz \right ],\label{Def-h(t)-0}\\
m(t)&:=\left[\frac{\p H_t}{\p t}(G_t) \pz \right ]. \label{Def-m(t)-0}
\end{align}

The first three components are given 
by the following proposition which was
proved in Section $3.2$ in \cite{GTS-II}. 

\begin{propo}\label{TaylorExpansion-DLog}
There exist unique $f(t), g(t), d(t) \in \cDazz[[t]]$ 
with $f(0)=1$ and $d(0)=0$
such that, for any $u_t(z)\in \kttzz$, we have
\begin{align}
f(-t)\, u_t(z)&=u_t(F_t),\label{NewTaylorExpansion-e1}\\
g(t)\, u_t(z)&=u_t(G_t), \label{NewTaylorExpansion-e2}\\
e^{d(t)}\, u_t(z) &=u_t(G_t), \label{NewDLog-e2}
\end{align}
where, as usual, the exponential in 
Eq.\,$(\ref{NewDLog-e2})$ is given by 
\begin{align}\label{L2.2.3-e2}
e^{d(t)} = \sum_{m\geq 0} \frac {d(t)^m}{m!}.
\end{align}
\end{propo}

Note that, when we write $d(t)$ above as 
$d(t)= -\lb a_t(z)\pz \rb$ for some 
$a_t(z)\in t\kttzz$, then we get the so-called 
{\it D-Log}\label{D-Log} $a_t(z)$ of the automorphism 
$F_t(z)\in \ataz$, which has been studied 
in \cite{E1}-\cite{E3}, \cite{Z-exp} and \cite{WZ} 
for the commutative case. 

\begin{theo} \label{S-Correspondence-2} $($\cite{GTS-II}$)$
For any $\alpha\geq 1 $ and $F_t(z)\in \ataz$, 
we have,

$(a)$ the $5$-tuple $\Omega_{F_t}$ defined 
in Eq.\,$(\ref{Def-Omega-Ft})$
forms a \cNcs  system over 
the $K$-algebra $\cDazz$.

$(b)$ there exists a unique homomorphism 
$\cS_{F_t}: \cNsf \to \cDazz$ of $K$-Hopf algebras such that 
$\cS_{F_t}^{\times 5}(\Pi)=\Omega_{F_t}$. 
\end{theo}

\subsection{A \cNcs  System over the Grossman-Larson Hopf Algebra 
of Labeled Rooted Trees} \label{S2.3}

In this subsection, we recall the 
\cNcs  system $(\cH_{GL}^W, \Omega_\bT^W)$
constructed in \cite{GTS-IV} over 
the Grossman-Larson Hopf algebra $\mathcal H^W_{GL}$ 
of rooted trees labeled by positive integers 
of an non-empty 
$W\subseteq \bN^+$.
First, let us fix the following notation 
which will be used throughout the rest of this paper.

\vskip2mm

{\bf Notation:}

\vskip2mm

By a {\it rooted tree} we mean a finite
1-connected graph with one vertex designated as its {\it root.} 
For convenience, we also view the empty set $\emptyset$ 
as a rooted tree and call it the {\it emptyset} rooted tree.
The rooted tree with a single vertex 
is called the {\it singleton} 
and denoted by $\circ$.
There are natural ancestral relations 
between vertices.  
We say a vertex $w$ is a {\it child} 
of vertex $v$ 
if the two are connected by an
edge and $w$ lies further from the root than $v$.  
In the same situation, we say $v$ 
is the {\it parent} of $w$.  
A vertex is called a {\it leaf}\/ if it has no
children. 

Let $W\subseteq \bN^+$ be 
a non-empty subset of positive integers. 
A {\it $W$-labeled rooted tree}
is a rooted tree 
with each vertex labeled by 
an element of $W$. If an element $m \in W$ 
is assigned to a vertex $v$, 
then $m$ is called the {\it weight} 
of the vertex $v$. When we speak of isomorphisms 
between unlabeled (resp.\,\,$W$-labeled) rooted trees, 
we will always mean isomorphisms 
which also preserve 
the root (resp.\,\,the root and also the labels of vertices).
We will denote by $\mathbb T$ (resp.\,\,$\mathbb T^W$) 
the set of isomorphism classes 
of all unlabeled (resp.\,\,$W$-labeled) rooted trees. 
A disjoint union of any finitely many rooted trees 
(resp.\,\,$W$-labeled rooted trees) 
is called a {\it rooted forest} (resp.\,\,$W$-labeled {\it rooted forest}).
We denote by $\mathbb F$ (resp.\,\,$\mathbb F^W$) 
the set of unlabeled (resp.\,\,$W$-labeled) rooted forests.  

\vskip2mm

With these notions in mind, we establish the following notation.

\begin{enumerate}
\item For any rooted tree $T\in \bT^W$, we set the following notation:
\begin{itemize}
\item $\text{rt}_T$ denotes the root vertex of $T$ and $O(T)$
the set of all the children of $\text{rt}_T$. We set
$o(T)=|O(T)|$ (the cardinal number of the set $O(T)$). 

\item $E(T)$ denotes the set of edges of $T$.

\item $V(T)$ denotes the set of  vertices of $T$ and $v(T)=|V(T)|$.

\item $L(T)$ denotes the set of leaves of $T$ and $l(T)=|L(T)|$

\item For any $v\in V(T)$ define the {\it height} of $v$ 
to be the number of edges in the (unique) geodesic connecting 
$v$ to $\text{rt}_{T}$. The {\it height} of $T$ is defined 
to be the maximum of the heights of
its vertices.

\item For any $T\in \bT^W$ and $T\neq \emptyset$, 
$|T|$ denotes the sum of the weights of 
all vertices of $T$. 
When $T=\emptyset$, we set $|T|=0$.

\item For any $T\in \bT^W$, we denote 
by $\text{Aut}(T)$ the automorphism group 
of $T$ and $\alpha(T)$ the cardinal 
number of $\text{Aut}(T)$.

\end{itemize}

\vskip2mm

\item Any subset of $E(T)$ is called a {\it cut} of $T$. 
A cut $C\subseteq E(T)$ is said to be {\it admissible} 
if no two different edges of $C$ lie in 
the path connecting the root and a leaf. We denote by 
$\mathcal C(T)$ the set of all admissible cuts of $T$. 
Note that, the empty subset $\emptyset$  of $E(T)$
and $C=\{e\}$  
for any $e\in E(T)$ 
are always admissible cuts. We will identify any edge 
$e\in E(T)$ with the admissible cut $C:=\{e\}$ 
and simply say the edge $e$ itself is 
an admissible cut of $T$.

\vskip2mm

\item For any 
$T \in \bT^W$ with $T\neq \circ$, 
let $C\in \mathcal C(T)$ 
be an admissible cut of $T$ 
with $|C|=m\geq 1$.
Note that, after deleting 
the edges in $C$ from $T$, 
we get a disjoint union of $m+1$ 
rooted trees, 
say $T_0$, $T_1$, ..., $T_m$ 
with $\text{rt}(T)\in V(T_0)$.
We define $R_C(T)=T_0 \in \mathbb T^W$ 
and $P_C (T)\in \mathbb F^W$ 
the rooted forest formed by $T_1$, ..., $T_m$. 

\vskip2mm

\item For any disjoint admissible cuts 
$C_1$ and $C_2$, we say
``$C_1$ lies above $C_2$", and write $C_1\succ
C_2$, if $C_2 \subseteq E(R_{C_1}(T))$. 
This merely says that all edges of $C_2$ remain 
when we remove all edges of $C_1$ and 
$P_{C_1}(T)$.  Note that this relation
is not transitive. 
When we write $C_1\succ \cdots \succ
C_r\,$ for $C_1,\ldots, C_r\in \mathcal C(T)$,   
we will mean that $C_i \succ C_j$ 
whenever $i<j$. 

\vskip2mm

\item For any $T\in \bT^W$, we say $T$ 
is a {\it chain} if 
its underlying rooted tree
is a rooted tree with a single leaf.
We say $T$ is a {\it shrub} if 
its underlying rooted tree
is a rooted tree 
of height $1$. 
We say $T$ is {\it primitive} if 
its root has only one child. 
For any $m\geq 1$, we set $\mathbb H_m$, 
$\mathbb S_m$ and $\mathbb P_m$ to be the sets of 
the chains, shrubs and primitive rooted trees 
$T$ of weight $|T|=m$, respectively. 
$\mathbb H$, $\mathbb S$ and $\mathbb P$ 
are set to be the unions of $\mathbb H_m$, 
$\mathbb S_m$ and $\mathbb P_m$, 
respectively, for all $m\geq 1$.   
\end{enumerate}

\vskip2mm

Let $K$ be any unital commutative $\bQ$-algebra 
and $W$ a non-empty subset of positive integers.
First, let us recall the Connes-Kreimer Hopf algebras 
$\mathcal H_{CK}^W$ of labeled rooted forests. 

As a $K$-algebra, the Connes-Kreimer Hopf algebra 
$\mathcal H_{CK}^W$ is the free commutative algebra 
generated by formal variables 
$\{X_T \,|\, T\in \mathbb T^W \}$. 
Here, for convenience, we will still use $T$ 
to denote the variable $X_T$ in 
$\mathcal H_{CK}^W$.
The $K$-algebra product is given by the disjoint union. 
The identity element of this algebra, denoted by $1$,  
is the free variable $X_\emptyset$ 
corresponding to the emptyset rooted tree. 
The coproduct $\Delta: \mathcal H_{CK}^W \to 
\mathcal H_{CK}^W \otimes \mathcal H_{CK}^W$ 
is uniquely determined by setting 
\begin{align}
\Delta(1)&=1\otimes 1, \label{CK-Delta-1} \\
\Delta(T)&= T\otimes 1+ 
\sum_{C\in \mathcal C(T)} P_C(T) \otimes R_C(T). \label{CK-Delta-2}
\end{align}
The co-unit $\epsilon: \mathcal H_{CK}^W \to K$ 
is the $K$-algebra homomorphism 
which sends $1\in \mathcal H_{CK}^W$ to $1\in K$ 
and $T$ to $0$ for any $T\in \mathbb T^W$ 
with $T\neq \emptyset$. 
With the operations defined above and 
the grading given by the weight, 
the vector space $\mathcal H_{CK}^W$ 
forms a graded commutative 
bi-algebra, hence there is a unique antipode 
$S: \mathcal H_{CK}^W\to \mathcal H_{CK}^W$ 
that makes $\mathcal H_{CK}^W$ a Hopf algebra.

Next we recall the Grossman-Larson Hopf algebra 
of labeled rooted trees. 
First we need define the following operations 
for labeled rooted forests.  For any 
labeled rooted forest $F$ which is disjoint 
union of labeled 
rooted trees $T_1$, $T_2$, ... , $T_m$, 
we set $B_+(T_1, T_2, \cdots, T_m)$ 
the rooted tree obtained by connecting roots 
of $T_i$ $(1\leq i\leq m)$ 
to a newly added root. We will keep the labels 
for the vertices of $B_+(T_1, T_2, \cdots, T_m)$ 
from $T_i$'s, but for the root, we label it by $0$. 

Furthermore, we fix the following convention 
for the operation $B_+$. 
First, for any $T_i\in \bT^W$  
$(1\leq i\leq m)$ and $j_i \geq 1$, the notation 
$B_+(T_1^{j_1}, T_2^{j_2}, \cdots, T_m^{j_m} )$
denotes the rooted tree obtained by applying 
the operation $B_+$ to $j_1$-copies of $T_1$; 
$j_2$-copies of $T_2$; and so on. Secondly, 
for any $m\geq 1$, we will extend the operation 
$B_+$ multi-linearly 
to a linear map $B_+$ from
$\lp \cH_{CK}^W \rp^{\times m}$ to the vector spaces
spanned by the resulted rooted trees. 
 
Now, we set $\bar {\mathbb T}^W:=\{ B_+(F) \, | \, F \in \mathbb F^W \}$.
Then,  $B_+: \mathbb F^W \to \bar{\mathbb T}^W$ becomes a bijection.
We denote by $B_- :  \bar{\mathbb T}^W \to \mathbb F^W$ 
the inverse map of $B_+$. More precisely, for any 
$T\in \bar{\mathbb T}^W$, 
$B_-(T)$ is the $W$-labeled rooted forest obtained by cutting off 
the root of $T$ as well as all edges connecting to the root in $T$.

Note that, precisely speaking,
elements of $\bar{\mathbb T}^W$ are not 
$W$-labeled trees for $0\not \in W$. 
But, if we set 
$\bar W=W\cup\{0\}$, 
then we can view $\bar{\mathbb T}^W$ 
as a subset 
of $\bar W$-labeled 
rooted trees $T$ with the root $\text{rt}_T$ 
labeled by $0$ 
and all other vertices 
labeled by non-zero elements of $\bar W$. 
We extend the definition of
the weight for elements 
of $\mathbb F^W$ to elements of 
$\bar{\mathbb T}^W$ by simply counting 
the weight of roots by zero. 
We set $\bar{\mathbb S}_m^W:=B_+(\mathbb S_m^W)$
$(m\geq 1)$ and $\bar{\mathbb S}^W:=B_+(\mathbb S^W)$. 
We also define
$\bar{\mathbb H}_m^W$, 
$\bar{\mathbb P}_m^W$, 
$\bar{\mathbb H}^W$ and 
$\bar{\mathbb P}^W$ 
in the similar way.

The Grossman-Larson Hopf algebra $\mathcal H_{GL}^W$ as a vector space 
is the vector space spanned by elements of $\bar{\mathbb T}^W$ 
over $K$. For any $T\in \bar{\mathbb T}^W$, we will still denote by 
$T$ the vector in $\mathcal H_{GL}^W$ that is corresponding to $T$. 
The algebra product is defined as follows. 
For any $T, S\in \bar{\mathbb T}^W$ with 
$T=B_+(T_1, T_2, \cdots, T_m)$, we set $T\cdot S$ to be the sum of  
the rooted trees obtained by connecting the roots of $T_i$ 
$(1\leq i\leq m)$ to vertices of $S$ 
in all possible $m^{v(S)}$ different ways. 
Note that, the identity element with respect to this 
algebra product is given by the singleton 
$\circ=B_+(\emptyset)$. But we will denote it by $1$.

To define the co-product $\Delta: \cH_{GL}^W \to \cH_{GL}^W \otimes \cH_{GL}^W$,
we first set 
\begin{align}
\Delta (\circ)=\circ \otimes \circ. 
\end{align}
 
Now let $T\in \bar{\mathbb T}^W$ with $T\neq \circ$, 
say $T=B_+ (T_1, T_2, \cdots, T_m)$ with
$m \geq 1$ and $T_i\in \mathbb T^W$ 
$(1\leq i\leq m)$. 
For any non-empty subset 
$I\subseteq \{1, 2, \cdots, m\}$,
we denote by 
$B_+(T_I)$ the rooted tree 
obtained by applying the $B_+$ operation 
to the rooted trees $T_i$ 
with $i\in I$. 
For convenience, when $I=\emptyset$, 
we set $B_+(T_I)=1$. 
With these notation fixed, the co-product 
for $T$ is given by
\begin{align}
\Delta (T)=\sum_{I\sqcup J=\{1, 2, \cdots, m\}} 
B_+(T_I) \otimes B_+(T_J).
\end{align}

Note that, a rooted tree in 
$\bar{\mathbb T}^W$ is a primitive element of 
the Hopf algebra $\cH_{GL}^W$ 
iff it is a primitive rooted tree 
in the sense that we defined before, 
namely the root of $T$ has 
one and only one child.

\begin{rmk}\label{action}
Note that, for any $S \in \bbT^W$ 
and $T \in \mathbb F^W$,
we also can define a 
``product'' still denoted by $S\cdot T$ in the 
exact same way as we define the product 
of elements of $\bbT^W$. 
By the linear extention, 
this ``product" makes 
$\cH_{CK}^W$ a $K$-algebra module of 
$\cH_{GL}^W$.
\end{rmk}

The following results later will be very useful 
in our later arguments. 

\begin{theo}\label{Dual-MM}
$(a)$
The Hopf algebras $\cH_{GL}^W$ 
and $\cH_{CK}^W$ are graded dual to each other. 
The pairing is given by, 
for any $T\in \bar{\mathbb T}^W$ 
and $S\in \mathbb F^W$, 
\begin{align}\label{pairing}
<T, F>=\begin{cases} 0, & \text{ if } T \not \simeq B_+(F),\\
\alpha (T), &\text{ if } T \simeq B_+(F).
\end{cases}
\end{align}

$(b)$ $\cH_{GL}^W$ as a Hopf algebra is 
isomorphic to the universal enveloping algebra of 
the Lie algebra formed by its primitive elements, 
which are exactly linear combinations of 
the primitive rooted trees. In particular, 
$\cH_{GL}^W$ as an $K$-algebra 
is generated by the 
primitive rooted trees. 
\end{theo}

For a proof of $(a)$, see \cite{H} and \cite{F2}.
$(b)$ follows directly from the well-known 
Milnor-Moore's Theorem (\cite{MM}), since $\cH_{GL}^W$ 
is a connected graded and cocommutative Hopf algebra.

Next, let us recall the following lemma proved 
in \cite{GTS-IV} that will be crucial for our later
arguments.

Let $\vec{C}=(C_1, \ldots, C_r)\in \mathcal C(T)^{\times r}$ 
be a sequence of admissible cuts
with $C_1 \succ \cdots \succ C_r$.
We define a sequence of $T_{\vec{C},1}, \ldots , 
T_{\vec{C},r+1} \in \bar{\mathbb T}^W$ 
as follows: we first set  $T_{\vec{C},1}= B_+ (P_{C_1}(T))$ and let
$S_1= R_{C_1}(T)$.  Note that 
$C_2, \ldots, C_r\in \mathcal C (S_1)$. 
We then set $T_{\vec C, 2}=B_+ (P_{C_2}(S_1))$ 
and  $S_2= R_{C_2}(S_1)$ and repeat this procedure 
until we get $S_{r}=R_{C_r}(S_{r-1})$ and then set 
$T_{\vec C, r+1}=S_r$.  
In the case that, each $C_i$ $(1\leq i\leq r)$
consists of a single edge, 
say $e_i\in E(T)$,  we simply denote $T_{\vec{C},i}$ 
by $T_{e_i}$.

\begin{lemma}\label{key-lemma-2}  
For any $r\geq 1$, 
$y=\{ y_T^{(i)}\,|\, 1\leq i\leq r; \,, T\in \bar{\bT}^W\}$
be a collection of commutative formal variables.
Then, we have, 
\allowdisplaybreaks{
\BQ\label{key-lemma2-e1}
&{}&
\sum_{\substack{(T_1,\ldots,T_r)\in (\bar{\bT}^W)^r }}
\left[ y_{T_1}^{(1)}\mathcal V_{T_1}\right]
\cdots \left[y_{T_r}^{(r)}\mathcal V_{T_{r}}\right] \\
&{}&\qquad\qquad \qquad =
\sum_{T\in \bar{\bT}^W }\,\,\,\,
\sum_{\substack{\vec{C}=(C_1,\ldots,C_r) \in {\mathcal C} (T)^r \\C_1\succ \cdots
\succ C_r}}y_{T_{\vec{C},1}}^{(1)}\cdots
y_{T_{\vec{C},r}}^{(r)}  \mathcal V_T. \nno 
\EQ}
\end{lemma}

Now we recall the \cNcs  system $\Omega_{\bT}^W$ constructed in 
\cite{GTS-IV} over the Grossman-Larson Hopf algebra $\cH_{GL}^W$.

First, let us define the following constants for the rooted trees in $\bbT^W$:
\begin{enumerate}
\item[$\bullet$] We set $\beta_T$ to be the weight 
of the unique leaf of $T$ if $T\in \bar{\mathbb H}^W$ and $0$ otherwise. 
\item[$\bullet$] We set $\gamma_T$ to be the weight 
of the unique child of the root of $T$ if $T\in \bar{\mathbb P}^W$ and $0$ otherwise. 
\item[$\bullet$]\label{theta} We set $\theta_T$ to be the coefficient of $s$ 
of the order polynomial $\Omega(B_-(T), s)$ of the underlying 
unlabeled rooted forest of $B_-(T)$.
\end{enumerate}

For general studies on the order polynomials $\Omega(P, s)$ 
of finite posets $P$, see \cite{St1}. 
For an interpretation of the constant
$\phi_T=(-1)^{v(T)-1}\varphi_T=(-1)^{v(T)-1}\theta_{B_+(T)}$  
in terms of the numbers of chains with fixed lengths in  
the lattice of the ideals of the poset $T$, 
see Lemma $2.8$ in \cite{SWZ}.

The following results\footnote{Note that, the approach 
to $\theta_T$'s in \cite{GTS-IV} is different from the one 
we adapt here. But it was shown there that 
the constants determined by the properties in 
Proposition \ref{propo-theta} 
are same as the $\theta_T$'s we defined here.}
on $\theta_T$ $(T\in \bbT^W)$ proved in Section 
$5$ in \cite{GTS-IV} will be needed later.

\begin{propo}\label{propo-theta}
\begin{enumerate}
\item For the singleton $\circ$ and 
any non-primitive rooted tree $T\in \mathbb T$, 
i.e. $o(T)>1$, we set $\theta_\circ = \theta_T=0$.

\item  For $T=B_+(\circ)$, we set $\theta_T=1$.

\item  For any primitive $T\in \mathbb P$ with $v(T)\geq 3$, we define $\theta_T$ 
inductively by 
\begin{align}\label{theta-Recur}
\theta_T = 1- \sum_{m \geq 2} \frac 1{m!} 
\sum_{\substack{\vec e=(e_1,\ldots,e_{m-1})\in E(T)^{m-1}\\e_1\succ
\cdots\succ e_{m-1}}}
\theta_{T_{e_1}} \theta_{T_{e_2}} \cdots  \theta_{T_{e_{m}}}. 
\end{align}
\end{enumerate}
\end{propo}

We will also need the following proposition proved
in \cite{GTS-IV} for the constants 
$\theta_T$ $(T\in \bT)$.

\begin{propo}\label{theta-orderP}
For any $T\in \mathbb P$, we have
\begin{align}\label{theta-orderP-e1}
\nabla \Omega(T, s) & = \theta_T \,  s +  \sum_{k=2}^{v(T)} \frac {s^{k}}{k!}
\sum_{\substack{\vec e=(e_1,\ldots,e_{k-1})\in E(T)^{k-1}\\e_1\succ
\cdots\succ e_{k-1}}}
\theta_{B_-(T_{\vec e,1})}  \cdots \theta_{B_-(T_{\vec e,k-1})} \theta_{T_{\vec e,k}},
\end{align}
where $\nabla: K[s]\to K[s]$ is the linear operator 
that maps any $f(s)\in K[s]$ to $f(s)-f(s-1)$.  
\end{propo}

Now we consider the following generating functions of $T\in \bbT^W$.

\allowdisplaybreaks{
\begin{align}
\tilde f(t):&=\sum_{T \in \bar{\mathbb S}^W} (-1)^{o(T) }  t^{|T|} {\mathcal  V_T}
=1 + \sum_{\substack{T \in \bar{\mathbb S}^W \\
T \neq \circ }} (-1)^{o(T) }  t^{|T|} {\mathcal  V_T}, 
\label{Def-3f(t)} \\ 
\tilde g(t):&= \sum_{T \in \bar{\mathbb T}^W} t^{|T|} {\mathcal  V_T}
=1 + \sum_{\substack{ T  \in \bar{\mathbb T}^W \\
T \neq \circ }} t^{|T|} {\mathcal  V_T}, 
\label{Def-3g(t)}\\ 
\tilde d(t):&= \sum_{T \in \bar{\mathbb P}^W}  t^{|T|} \theta_T {\mathcal  V_T}. 
\label{Def-3d(t)}\\
\tilde h(t):&= \sum_{T \in \bar{\mathbb H}^W} t^{|T|-1}  \beta_T \mathcal  V_T, 
\label{Def-3h(t)} \\ 
\wtilde m(t):&= \sum_{T \in \bar{\mathbb P}^W}  t^{|T|-1} \gamma_T {\mathcal  V_T}, 
\label{Def-3m(t)}
\end{align} 
where, for any $T\in \bbT^W$, $\mathcal V_T:=\frac 1{\alpha(T)} T.$
}
We further set
\begin{align}\label{Def-3-Omega}
\Omega_\bT^W:=
(\, \tilde f(t),\, \tilde g(t), 
\, \tilde d\,(t), \, \tilde h(t), \wtilde m(t)\, ). 
\end{align}

\begin{theo}\label{Main-Thm-Trees} $($\cite{GTS-IV}$)$\,
 For any non-empty set $W\subseteq\bN^+$, we have

$(a)$ the $5$-tuple $\tilde \Omega_{F_t}$ defined 
in Eq.\,$(\ref{Def-3-Omega})$
forms a \cNcs  system over 
the Grossman-Larson Hopf algebra $\mathcal H^W_{GL}$.

$(b)$ there exists a unique homomorphism 
$\cT_W: \cNsf \to \mathcal H^W_{GL}$ of 
graded $K$-Hopf algebras such that 
$\cT^{\times 5}_W(\Pi)=\tilde \Omega_{F_t}$. 
\end{theo}

\renewcommand{\theequation}{\thesection.\arabic{equation}}
\renewcommand{\therema}{\thesection.\arabic{rema}}
\setcounter{equation}{0}
\setcounter{rema}{0}

\section{\bf A Hopf Algebra 
Homomorphisms from $\cH_{GL}^W$ to $\cDazz$}  
\label{S3}

In this section, for any non-empty $W\subseteq\bN^+$, 
$\alpha \geq 1$ and $F_t\in \ataz$ 
satisfying Eq.\,(\ref{Def-Hm}) below,
we construct a $K$-Hopf algebra 
homomorphism $\mathcal A: \cH_{GL}^W\to \cDazz$ such that
$\mathcal A^{\times 5}$ maps the \cNcs  system $\Omega_\bT^W$ 
in Theorem \ref{Main-Thm-Trees} to the \cNcs  system $\Omega_{F_t}$ 
in Theorem \ref{S-Correspondence-2}.

Let $K$, $z$ and $t$ be as given 
in Subsection \ref{S2.2}. We will also 
use the notation fixed in 
Sections \ref{S3} freely throughout 
this section.

Let us start with the introducing of 
the following two operations for 
the $K$-derivations of $\kzz$.

First, for any $\phi, \delta \in \cDrzz$ with 
$\delta=\lb f(z)\pz \rb$, we set
\begin{align}\label{collapse}
\phi \triangleright \delta:=\lb (\phi f)(z)\pz\rb. 
\end{align}

To define the second operation, let
$w=(w_1, w_2, \cdots, w_n)$ be 
another $n$ free variables which are 
independent with the free variables 
$z$.  For any $K$-derivations 
$\delta_i=\lb \Vec v_i(z)\pz \rb$ 
$(1\leq i\leq m)$ with $\Vec v_i(z)\in \kzz^{\times n}$, 
we define 
$B_+(\delta_1, \delta_2, \cdots, \delta_m)$ 
by setting, for any $u(z)\in \kzz$, 
\begin{align}\label{Def-B+}
& B_+(\delta_1, \delta_2, \cdots, \delta_m)u(z):=\\
&\quad\quad\quad\quad
\left. \lb \Vec v_1(w)\pz\rb  \lb \Vec v_2(w)\pz \rb\cdots 
\lb \Vec v_m(w)\pz \rb u(z)\right |_{w=z}.\nno 
\end{align}

Furthermore, for any $k_i\geq 0$ $(1\leq i\leq m)$, 
we let 
$B_+(\delta_1^{k_1}, \delta_2^{k_2}, \cdots, \delta_m^{k_m})$ 
denote the operator obtained by applying $B_+$ 
to the multi-set of $j_1$-copies of $\delta_1$; 
$j_2$-copies of $\delta_2$, ..., 
$j_m$-copies of $\delta_m$.

Note that $B_+(\delta_1, \delta_2, \cdots, \delta_m)$ 
is multi-linear and symmetric on $\delta_i$ 
$(1\leq i\leq m)$. 
When $m=1$, $B_+(\delta_1)=\delta_1$.

The following two lemmas have been proved in 
Section $3$ in \cite{GTS-II}.

\begin{lemma}\label{LL4.1.1}
$(a)$ Let $\delta_i \in \cDrzz$ $(1\leq i\leq m)$. 
Then, for any $\phi\in \cDzz$, we have
\begin{align}\label{LL4.1.1-e1}
\phi \cdot B_+(\delta_1, \delta_2, \cdots, \delta_m)&=
B_+(\phi, \delta_1, \delta_2, \cdots, \delta_m) \\ 
& \quad +\sum_{i=1}^m B_+
(\delta_1, \cdots, \phi\triangleright \delta_i, \cdots  \delta_m).\nno
\end{align}

$(b)$ For any $\delta_i \in \cDrzz$ $(1\leq i\leq m)$, 
$B_+(\delta_1, \delta_2, \cdots, \delta_m) \in \cDzz$. 
\end{lemma}

\begin{lemma}\label{Taylor-f(t)}
In terms of the $B_+$ operation defined above, 
$f(t)\in \cDazz$ defined in Proposition 
$(\ref{TaylorExpansion-DLog})$ is given by
\begin{align}\label{Taylor-f(t)-e1}
f(t)= \sum_{k \geq 0} \frac {(-1)^{k}}{k!}
B_+\lp \lb H_t(z) \pz \rb^k \rp.
\end{align}
\end{lemma}

Now, we fix a non-empty subset $W \subseteq \bN^+$ and
$\alpha \geq 1$. Let $F_t(z)=z-H_t(z) \in \ataz$ 
such that 
\begin{align}\label{Def-Hm}
H_t(z)=\sum_{m \in W} t^m H_{[m]}(z),
\end{align}
for some $H_{[m]}(z)\in \kzz^{\times n}$ $(m\in W)$.


First, we assign $P_T(z)\in \kzz^{\times n}$ for each 
$T\in \bT^W$ inductively as follows.
\begin{enumerate}

\item[$(1)$] For the singleton $\circ$ labeled by $m \in W$, 
denoted by $\circ_m$, we set $P_{\circ_m}(z):= H_{[m]}(z)$.

\item[$(2)$] For any non-singleton $T\in \mathbb T^W$ with $\text{rt}_T$ 
labeled by $m\in W$, write $T=B_+(T_1, T_2, ... , T_d)$ 
with $T_i\in \bT^W$ $(1\leq i\leq d)$ and 
\begin{align}\label{Def-PT}
P_T (z):= B_+ \lp \lb P_{T_1}(z)\pz\rb, 
... , \lb P_{T_m}(z) \pz \rb \rp H_{[m]}(z).
\end{align}
\end{enumerate}

Let $\cH_{CK}^W[1]$ be the vector subspace of $\cH_{CK}^W$ 
spanned by $T\in \bT^W$. We define a linear map 
$\mathcal U_{F_t}: \cH_{C}^W[1] \to \cDazz$
by setting
\begin{align}\label{Def-cU}
\mathcal U_{F_t}: \cH_{CK}^W[1] & \to \kzz^{\times n} \\
  T& \to P_T(z).\nno
\end{align}

Next we assign $D_T(z)\in \cDazz$ for each 
$T\in \bbT^W=B_+(\bT^W)$ as follows.

\begin{enumerate}
\item[$(1)$] For the singleton $T=\circ$, we set $D_T=id$. 

\item[$(2)$] For any rooted tree $T=B_+(T_1, T_2, \cdots, T_m)$ 
with $T_i\in \bT^W$, we set
\begin{align}\label{Def-DT}
D_T(z) :=  B_+ \lp \lb P_{T_1}(z) \pz\rb, \lb P_{T_2}(z) \pz \rb,
 ... , \lb P_{T_m}(z) \pz \rb \rp.
\end{align} 
\end{enumerate}

Note that, from the definition of $P_T(z)$'s above, 
it is easy to see inductively that, 
for any $T\in \bT^W$, $o(P_T(z))\geq \alpha$. 
Hence we do have $D_T \in \cDazz$ for any $T\in \bar{\bT}^W$.

Now we define a linear map $\mathcal A_{F_t}: \cH_{GL}^W \to \cDazz$
by setting
\begin{align}\label{Def-cA}
\mathcal A_{F_t}: \cH_{GL}^W & \to \cDazz \\
  T&\to D_T(z).\nno
\end{align}

When $F_t\in \ataz$ is clear in the context, 
$\mathcal U_{F_t}$ and $\cA_{F_t}$ will also be simply 
written as $\mathcal U$ and $\cA$, respectively. 

From the definitions above, the following 
lemma follows immediately.

\begin{lemma}\label{LL4.1.2}
$(a)$ For any $T\in \bT^W$ with the root 
labeled by $m\in W$, we have
\begin{align} \label{LL4.1.2-e1}
\cA (B_+(T))H_{[m]}(z) = D_T. 
\end{align}

$(b)$ For any $m\geq 1$,  
we have
\begin{align} \label{LL4.1.2-e2}
\cA \circ B_+ = B_+ \circ \cA^{\times m} \circ B_+^{\times m}, 
\end{align}
as maps from $(\bT^W)^{\times m}$ to $\cDazz$. 
\end{lemma}

Now let us prove the following technic lemma.

\begin{lemma}\label{LL4.1.3}
For any $S \in \bar{\mathbb P}^W$ and $T\in \bT^W$, we have
\begin{align}\label{LL4.1.3-e1}
D_S \cdot P_{T}= \mathcal U (S \cdot T),
\end{align}
where $S \cdot T$ is the ``product" of $S$ 
and $T$ as in Remark $(\ref{action})$.
\end{lemma}

\pf We use the induction on the number $v(T)$
of vertices of $T$. 

First, when $v(T)=1$, then $T$ is the singleton labeled 
by some $m\in W$. Then, by Eqs.\,$(\ref{Def-PT})$ and 
$(\ref{Def-DT})$, it is easy to see that
both sides of Eq.\,$(\ref{LL4.1.2-e2})$ 
in this case are $P_S(z)$.

Now, assume $v(T)\geq 2$. We write $S=B_+(S')$
and $T=(T_1, \cdots, T_d)$ with $S, T_i\in \bT^W$ 
$(1 \leq i \leq d)$. 
Let $m\in W$ be the label of the root of $T$.
Set $\phi:=D_S$ and $\delta_i:=D_{T_i}$ $(1\leq i\leq d)$. 
Note that, by our conditions, 
$\phi$ and  $\delta_i=D_{T_i}$ $(1\leq i\leq d)$
are all $K$-derivations. 
Now by Eqs.\,(\ref{Def-PT})
and (\ref{LL4.1.1-e1}), we have
\begin{align*}
D_S \cdot P_{T}&=
\phi \cdot B_+ ( \delta_1, \, \delta_2, \cdots \, \delta_d )\,  H_{[m]}(z) \\
&=B_+( \phi, \, \delta_1, \, \delta_2,\, \cdots \, \delta_d  )\, H_{[m]}(z) \\
& \quad \quad + \sum_{1\leq i\leq d} B_+( \delta_1, \cdots, \, 
\phi \triangleright \delta_i, \cdots, \, \delta_d  ) \, H_{[m]}(z) \\
\intertext{Applying Eq.\,(\ref{collapse}) and the induction assumption :}
&=B_+( \phi, \, \delta_1, \, \delta_2,\, \cdots \, \delta_d  )\, H_{[m]}(z) \\
& \quad \quad + \sum_{1\leq i\leq d} B_+( \delta_1, \cdots, \, 
\lb P_{S'\cdot T_i}(z)\pz \rb, \cdots, \, \delta_d  ) \, H_{[m]}(z) \\
\intertext{Applying Eq.\,(\ref{LL4.1.2-e2}):}
&=\cA \lp  B_+( S', T_1, \, \cdots, T_2 ) \rp  H_{[m]}(z)\\
& \quad  + \sum_{1\leq i\leq d} \cA
\lp \, B_+( T_1, \cdots, \, 
S'\cdot T_i, \cdots, \, T_d  \rp \, H_{[m]}(z) \\
\intertext{Applying Eq.\,(\ref{LL4.1.2-e1}):}
&=\mathcal U \lp B_+( S', T_1, \, \cdots, T_2 ) \rp \\
&\quad\quad + 
 \sum_{1\leq i\leq d}\mathcal U \lp B_+( T_1, \cdots, \, 
S'\cdot T_i, \cdots, \, T_d  )\rp
\end{align*}

Note that, by the definition of $S\cdot T$ (see Remark \ref{action}), 
we have
\begin{align}\label{LL4.1.3-pe1}
S\cdot T&= S\cdot B_+(T_1, \,T_2, \, \cdots, \, T_2 ) \\
& = B_+( S', T_1, \, \cdots, T_2 )+ \sum_{1\leq i\leq d}  
B_+( T_1, \cdots, \, S'\cdot T_i, \cdots, \, T_d  ). \nno
\end{align}

Combining the two equations above and using the linearity of $\mathcal U$, 
we get Eq.\,(\ref{LL4.1.3-e1}).
\epfv

Now, we can formulate and prove 
the first main result of this section.

\begin{theo}\label{S3-Main-1}
The linear map $\mathcal A_{F_t}: \cH_{GL}^W \to \cDazz$ 
is a homomorphisms of $K$-Hopf algebras.
\end{theo}

\pf Let us first show the linear map $\mathcal A$ 
is a homomorphism of $K$-algebras. By the definition 
of $D_T$'s and $\cA$, $\cA$ maps the identity element 
$\circ \in \cH_{GL}^W$ to the identity element 
$D_\circ=1$ of $\cDazz$. So we only need show, 
for any $S, T\in \bbT^W$, we have
\begin{align}\label{S3-Main-pe1}
D_S\cdot D_{T}=D_{S \cdot T},
\end{align}
where $\cdot$ denotes the both algebra product of 
$\cH_{GL}^W$ and $\cDazz$.

By the fact pointed in $(b)$ of Theorem \ref{Dual-MM}, 
we may assume that $S$ is primitive, in which case $D_S$ 
is a derivation.
Now, assume $v(T)\geq 2$. (When $v(T)=1$, Eq.(\ref{S3-Main-pe1}) 
is trivial).
We write $S=B_+(S')$ and $T=(T_1, \cdots, T_d)$ with $S, T_i\in \bT^W$ 
$(1 \leq i \leq d)$. Let $m\in W$ be the label of the root of $T$.
Set $\phi:=D_S$ and $\delta_i:=D_{T_i}$ $(1\leq i\leq d)$. 
Now by Eqs.\,(\ref{Def-DT})
and (\ref{LL4.1.1-e1}), we have
\begin{align*}
& D_S \cdot D_{T}=
\phi \cdot B_+ ( \delta_1, \, \delta_2, \cdots \, \delta_d ) \\
&   =B_+( \phi, \, \delta_1, \, \delta_2,\, \cdots \, \delta_d  ) 
  + \sum_{1\leq i\leq d} B_+( \delta_1, \cdots, \, 
\phi \triangleright \delta_i, \cdots, \, \delta_d  )  \\
\intertext{Applying Eqs.\,(\ref{collapse}) and (\ref{LL4.1.3-e1}):}
&=B_+( \phi, \, \delta_1, \, \delta_2,\, \cdots \, \delta_d  )
+ \sum_{1\leq i\leq d} B_+( \delta_1, \cdots, \, 
\lb P_{S'\cdot T_i}(z)\pz \rb, \cdots, \, \delta_d  )  \\
\intertext{Applying Eq.\,(\ref{LL4.1.2-e2}) and the linearity of $\cA$:}
&=\cA \lp  B_+( S', T_1, \, \cdots, T_2 ) 
 + \sum_{1\leq i\leq d}   B_+( T_1, \cdots, \, 
S'\cdot T_i, \cdots, \, T_d  \rp \, \\
\intertext{Applying Eq.\,(\ref{LL4.1.3-pe1}):}
&=\cA (S\cdot T) \\
&=D_{S\cdot T}.
\end{align*}

Hence, we have proved that $\cA$ is a $K$-algebra 
homomorphism. To see that it is also a homomorphism 
of Hopf algebras, note that, again, by the fact in $(b)$ of 
Theorem \ref{Dual-MM}, we know that $\cH_{GL}^W$ 
as a Hopf algebra is isomorphic to the universal enveloping algebra 
of the Lie sub-algebra $\cP(\cH_{GL}^W)$ generated 
by the primitive rooted trees. 
On the other hand, the Hopf algebra 
structure of $\cDazz$ is also obtained by viewing it as  
the universal enveloping algebra of the Lie algebra $\cDrzz$. 
Since the Lie brackets of $\cP(\cH_{GL}^W)$ and $\cDrzz$ are 
both given respectively by the commutator brackets of 
the algebra products of $\cH_{GL}^W$ and $\cDazz$, 
the restriction of $\cA: \cP(\cH_{GL}^W)\to \cDrzz$ is a 
homomorphism of Lie algebras. 
Consequently, $\cA$ also preserves the Hopf algebras of 
the corresponding enveloping algebras, 
which are $\cH_{GL}^W$ and $\cDazz$.
\epfv

Our second main result of this section is the following theorem.

\begin{theo}\label{S3-Main-2}
Let $\Oft$ and $\Omega_W$ the \cNcs systems in 
Theorems $\ref{S-Correspondence-2}$ and 
$\ref{Main-Thm-Trees}$, 
respectively, and  $\mathcal A_{F_t}: \cH_{GL}^W \to \cDazz$ 
be the $K$-Hopf algebras homomorphism 
in Theorem \ref{S3-Main-1}.
Then we have
$\mathcal A_{F_t}^{\times 5}(\Omega_\bT^W) = \Omega_{F_t}$. 
\end{theo}

Before we prove the theorem above, 
we need the following lemma proved 
in \cite{GTS-IV}, which gives a different way to look at 
the generating function $\tilde f(t)$ defined in 
Eq.\,(\ref{Def-3f(t)}).

\begin{lemma} \label{Taylor-for-3f(t)}
For any $m\in W$, let $\kappa_m$ denote the singleton 
labeled by $m$ and set
\begin{align}\label{Def-kappa}
\kappa(t):=\sum_{m\in W} t^m \kappa_m.
\end{align}
Then, we have
\begin{align}\label{Taylor-for-3f(t)-e1}
\tilde f(t)=1+\sum_{d\geq 1} \frac{(-1)^{d}}{d!} B_+( \, \kappa (t)^d \,),
\end{align}
where $B_+( \, \kappa (t)^d \,)$ denotes the term obtained 
by applying $B_+$ to $d$-copies of $\kappa (t)$.
\end{lemma}

\underline{Proof of Theorem \ref{S3-Main-2}}:
By Corollary $2.8$ in \cite{GTS-I}, 
it will be enough to show that $\cA$ maps 
one component of $\Omega_\bT^W$ to 
the component of $\Omega_{F_t}$ 
at the same location. Below we will 
show  $\cA(\tilde f(t))=f(t)$.


First, for any $m\in W$, let $\kappa_m$ denote 
the singleton labeled by $m$ and set
$\kappa(t):=\sum_{m\in W} t^m \kappa_m.$

By Eq.\,$(\ref{Def-Hm})$ and 
the definition of $\cA$ in Eq.\,$(\ref{Def-cA})$, 
we have
\begin{align}\label{S3-Main-pe2}
\cA B_+( \kappa(t)) )&=\sum_{m\in W} t^m \cA (B_+(\kappa (t)))\\
&=\sum_{m\in W} t^m \lb H_{[m]}(z)\frac \p{\p z} \rb  \nno \\
&=\lb H_t(z)\frac \p{\p z} \rb\nno
\end{align}

By Eq.\,$(\ref{LL4.1.2-e2})$ and the equation above, we have
for any $d\geq 1$, 
\begin{align}\label{S3-Main-pe3}
(\cA\circ B_+) ( \kappa(t))^d )&= (B_+ \circ \cA^d \circ B_+^d )(\kappa(t)^d)\\
&=B_+( \lb H_t(z)\frac \p{\p z} \rb^d ). \nno
\end{align}

Therefore, by Eq.\,$(\ref{Taylor-for-3f(t)-e1})$ in Lemma \ref{Taylor-for-3f(t)}
and the equation above, we have

\begin{align*}
\cA \lp \tilde f(t)\rp &=1+\sum_{d\geq 1} \frac{(-1)^{d}}{d!} \cA 
\lp B_+( \, \kappa (t)^d \,) \rp \\
&=1+\sum_{d\geq 1} \frac{(-1)^{d}}{d!} B_+( \lb H_t(z)\frac \p{\p z} \rb^d ) \\
\intertext{By Eq.\,$(\ref{Taylor-f(t)-e1})$ in Lemma \ref{Taylor-f(t)}:}
&=f(t). 
\end{align*}
\epfv

By applying Lemma \ref{LL4.1.3} and Theorem \ref{S3-Main-1}, it is easy 
to see the following proposition also holds.  

\begin{propo}\label{CD-1}
Let $\tilde {\mathcal L}: \cH_{GL}^W \times \cH_{CK}^W[1] \to \cH_{CK}^W[1]$
be the action induced by the natural action of elements $S\in \bbT^W$ 
on elements $T\in\bT^W$ $(\text{see Remark } \ref{action})$ and 
$\mathcal L: \cDazz \times \kzz \to \kzz^{\times n}$ 
the natural action of $\cDazz$ on $\kzz^{\times n}$. 
With $W$ and $F_t\in \ataz$ as before, 
 the following diagram commutes.
\begin{align}
\begin{CD}\label{CD-1-e1}
\cH_{GL}^W \times  \cH_{CK}^W[1] @>\tilde{\mathcal L} >> {\mathcal H}_{CK}^W[1] \\
@V \mathcal A_{F_t} \times \mathcal U_{F_t} VV @V  \mathcal U_{F_t}  VV\\
\cDazz \times \kzz^{\times n} @> \mathcal L >> { \quad \kzz^{\times n} }
\end{CD}
\end{align}
\end{propo}

Combining Theorems \ref{S-Correspondence-2}, \ref{Main-Thm-Trees} 
and \ref{S3-Main-1}, we have the following proposition.

\begin{propo}\label{CD-2}
For any $\alpha \geq 1$, let $W \subseteq\bN^+$ and 
$F_t\in \ataz$ fixed as before,  
we have the following commutative diagrams of 
$K$-Hopf algebra homomorphisms.
\begin{align}
\begin{CD}\label{CD-2-e1}
\cNsf @>{\cT_W} >> {\mathcal H}_{GL}^W \\
@V \cS_{F_t} VV @V \mathcal A_{F_t}  VV\\
\cDazz @= \cDazz
\end{CD}
\end{align}
\end{propo}

\pf By Theorems \ref{S-Correspondence-2}, \ref{Main-Thm-Trees} 
and \ref{S3-Main-2}, it is easy to see that 
\begin{align}
(\mathcal A_{F_t} \circ \cT_W)^{\times 5}\, (\Pi) =
\mathcal S_{F_t}^{\times 5}\, (\Pi) =\Omega_{F_t}. 
\end{align}

In particular, we have 
$\mathcal A_{F_t} \circ \cT_W (\lambda (t))= 
\mathcal S_{F_t} (\lambda (t))=f(t)$. 
Hence, for any $m\geq 1$, we have 
$\mathcal A_{F_t} \circ \cT_W (\Lambda_m)= 
\mathcal S_{F_t} (\Lambda_m)$. Since $\cNsf$ is the free 
$K$-algebra generated by $\Lambda_m$ $(m\geq 1)$, 
we have $\mathcal A_{F_t} \circ \cT_W = 
\mathcal S_{F_t}$.
\epfv

Next, let us consider the question when 
we can dualize the commutative diagram in the 
proposition above. First, we have to know when
the Hopf algebra 
homomorphism $\cS_{F_t}: \cNsf \to \cDazz$  
preserves the gradings of $\cS_{F_t}$ and $\cDazz$.
Note that, precisely speaking, $\cDazz$ 
is not graded in the usual sense, 
for some infinite sums are allowed in $\cDazz$. 
But we can consider the following graded 
subalgebras of $\cDazz$. 

Let $\cDz$ be the differential operator 
algebra of the polynomial
algebra $\kz$, i.e. $\cDz$ is the unital subalgebra of 
$\text{End}_K(\kz)$ generated 
by all $K$-derivations of $\kz$.  
For any $m\geq 0$, let $\cD_{[m]}\langle z \rangle$ 
be the set of all differential operators 
$U$ such that, for any homogeneous polynomial 
$h(z)\in \kz$ of degree $d\geq 0$, $U h(z)$ 
either is zero or is homogeneous of degree $m+d$. 
For any $\alpha\geq 1$, set $\cDaz:=\cDz\cap \cDazz$.
Then, we have the grading
\begin{align}\label{Grading-cDz}
\cDaz &=\bigoplus_{m \geq \alpha-1} 
\cD_{[m]} \langle z \rangle, 
\end{align}
with respect to which $\cDaz$ 
becomes a graded $K$-Hopf algebra.

Now, for any $\alpha \geq 2$, we let 
$\gtaz$ be the set of all automorphisms 
$F_t\in \ataz$ such that $F_t(z)=t^{-1}F(tz)$
for some automorphism $F(z)$ of $\kzz$.
It is easy to check that $\gtaz$ 
is a subgroup of $\ataz$. 
Then we have the following proposition 
proved in \cite{GTS-II}.

\begin{propo}\label{cS-graded}
For any $\alpha\geq 2$ and $F_t\in \ataz$, 
the differential operator specialization 
$\cS_{F_t}$ is a graded $K$-Hopf algebra 
homomorphism $\cS_{F_t}: \cNsf\to \cDaz\subset \cDazz$ 
iff $F_t\in \gtaz$.
\end{propo}

Now, for any $F_t\in \gtaz$ $(\alpha\geq 2)$, 
by the proposition above, 
we can take the graded dual of 
the graded $K$-Hopf algebra homomorphism 
$\cSft:\cNsf \to \cDaz$ 
and get the following corollary.

\begin{corol}\label{S-cQsf}
For any $\alpha \geq 2$ and $F_t\in\gtaz$, 
let $\cDaz^*$ be the graded dual 
of the graded $K$-Hopf algebra $\cDaz$. Then, 
$$
\cSft^*: \cDaz^* \to \cQf
$$ 
is a homomorphism of 
graded $K$-Hopf algebras. 
\end{corol}

By combining Proposition \ref{S-cQsf}
and Proposition \ref{CD-2} above, 
we have the following proposition.

\begin{propo}\label{CD-3}
For any $\alpha \geq 2$ and $F_t\in\gtaz$, 
we have the following commutative diagrams of graded 
$K$-Hopf algebra homomorphisms.
\begin{align}
\begin{CD}
{\mathcal Q}Sym @ < {\mathcal T^*} << {\mathcal H_{CK}^W } \\
@A {\cS_{F_t}^*} AA @A {\mathcal A^*_{F_t} } AA  \\
{\cDaz^* } @= {\cDaz^*}
\end{CD}
\end{align}
\end{propo}

\renewcommand{\theequation}{\thesection.\arabic{equation}}
\renewcommand{\therema}{\thesection.\arabic{rema}}
\setcounter{equation}{0}
\setcounter{rema}{0}

\section[{\bf Tree Expansion Formulas}]
{\bf Tree Expansion Formulas for D-Log's and Formal Flows}\label{S4}

In this section, we apply the $K$-Hopf algebra homomorphisms 
$\mathcal A_{F_t}: \cH_{GL}^W \to \cDazz$ ($F_t\in \ataz$) 
constructed in Section \ref{S3} to 
derive tree expansion formulas for the D-Log, 
the formal flow and the inverse map of $F_t$.
Note that, for the commutative case,
the tree expansion formula Eq.\,$(\ref{S4-Main-C-e3})$ 
for the inverse map was first given in \cite{BCW} 
and \cite{Wr3}. Later, it was generalized 
in \cite{WZ} to the D-Log's and the formal flows 
(see Eqs.\,$(\ref{S4-Main-C-e1})$ 
and $(\ref{S4-Main-C-e2})$). 
The proofs given here do not depend on 
the commutativity of the free variables. 
It not only generalizes 
the tree expansion formulas in 
\cite{BCW}, \cite{Wr3} and \cite{WZ} 
for the inverse maps,  the D-Log's and the formal flows 
to the noncommutative case, 
but also provides some 
new understandings to 
these formulas from 
the \cNcs  system point view.


First, we let $K$, $z$ and $t$ as before and 
fix an automorphism $F_t(z)\in \ataz$ 
$(\alpha \geq 1)$. As before, we fix the following notation 
for $F_t$ and its inverse map $G_t:=F_t^{-1}$:
\begin{align}
F_t(z)&=z-H_t(z),\\
G_t(z)&=z+tN_t(z),
\end{align}
with $H_t(z), N_t(z)\in \ktzz^{\times n}$.

Note that, in terms of the notation in Section \ref{S3}, 
we have $M_t(z)=tN_t(z)$.
Recall that, the D-Log of $F_t$ by definition 
is the unique $a_t(z)\in \kttzz$
such that, for any $u_t(z)\in \kttzz$, 
\begin{align}\label{Def-DLog}
e^{\lb a_t(z)\pz \rb }\cdot u_t(z)= u_t(F_t).
\end{align}

By the comments after Eq.\,(\ref{NewDLog-e2}), 
the relation of the D-Log $a_t(z)$ 
of $F_t$ with the third component $d(t)$ of 
the \cNcs  system $\Omega_{F_t}$ in Section \ref{S2.2} 
is given by 
\begin{align}
d(t)= -\lb a_t(z)\pz \rb. \label{d(t)-DLog}
\end{align}

Now let $s$ be another central parameter, i.e. 
it commutes with $z$ and $t$. We define
\begin{align}
F_t(z, s):= e^{s \lb a_t(z)\pz \rb }\,z= e^{-s\, d(t) }\,z. 
\label{Def-flow}
\end{align}

Note that, since $d(0)=0$, the exponential above 
is always well-defined. Actually it is easy to see 
$F_t(z, s)\in (K[s][[t]])
\langle \langle z \rangle \rangle^{\times n}$.
Therefore, for any $s_0 \in K$, $F_t(z, s_0)$ 
makes sense. 

Following its analog in \cite{E1}-\cite{E3} and 
\cite{WZ} in the commutative case, 
we call $F_t(z, s)$ the {\it formal flow} generated by $F_t(z)$.

Two remarks on the formal flow defined above are as follows.

First, it is well-known that 
the exponential of a derivation  
of any $K$-algebra $A$, when it makes sense, 
is always 
an automorphism of the algebra, so in our case, 
for any $s_0\in K$, $e^{s_0 \lb a_t(z)\pz \rb }$ 
is also an automorphism of $\kttzz$ over $K[[t]]$ 
which maps $z$ to $F_t(z, s_0)$. 
From Eq.\,$(\ref{Def-flow})$, 
it is clear that this automorphism also lies
in $\ataz$ since $o(a_t(z)) \geq \alpha$.

Secondly, by Eq.\,$(\ref{Def-flow})$ 
and the remark above, 
the formal flow $F_t(z, s)$ has 
the following properties:
\begin{align}
F_t(z, 0)& = z, \\
F_t(z, 1)&=F_t(z),\\
F_t( F_t(z, s_2), s_1)&=F_t(z, s_1+s_2),
\end{align}
for any $s_1, s_2\in K$.

In other words, $F_t(z, s)$ forms an one-parameter 
subgroup of the group $\ataz$. 
Therefore, for any integer $m\in \bZ$, 
$F_t(z, m)$ gives the $m^{th}$ 
(composing)\label{composingPower} power of $F_t$ 
as an element of the group $\ataz$. 
In particular, by setting $m=-1$, 
we get the inverse map $G_t$ of $F_t$, 
i.e. $F_t(z, -1)=G_t(z)$.

Now we consider the tree expansion formulas for 
the D-Log and formal flow of $F_t\in \ataz$. 
For convenience,  we first introduce the following short notations.
For any labeled rooted trees $T \in \mathbb T^W$ and $T' \in \bar{\mathbb T}^W$, 
we set 
\begin{align}
\cP_T(z)&=\mathcal U (\mathcal V_{T'})=\frac 1{\alpha (T)} P_T(z),\label{Def-cP} \\
\cD_{T'} (z)&= \mathcal A(\mathcal V_{T'})=\frac 1{\alpha (T')}D_{T'}(z).\label{Def-cD}
\end{align}

By Eqs.\,$(\ref{Def-PT})$ and $(\ref{Def-DT})$, 
it is easy to see that, 
for any primitive $T\in \bbT^W$, we have
\begin{align}\label{D-P}
\cD_{T} \cdot z=\cP_{B_-(T)}(z).
\end{align}

The main results of this subsection is the following theorem.

\begin{theo}\label{S4-Main}
For any $\alpha \geq 1$ and $F_t\in \ataz$, 
let $W$ be the set of positive integers such that 
$F_t$ can be written as in Eq.\,$(\ref{Def-Hm})$. Then

$(a)$ The D-Log of $F_t(z)$ is given by 
\begin{align}  \label{S4-Main-e1}
a_t(z)= \sum_{T\in {\mathbb T}^W} t^{|T|-1} \varphi_T {\mathcal P_T}(z), 
\end{align}
where, for any $T \in \mathbb T^W$, $\varphi_T$ is the 
coefficient of $s$ in the order polynomial 
$\Omega (T, s)$ of $T$ as an unlabeled rooted tree.

$(b)$ The formal flow $F_t(z, s)$ generated by $F_t(z)$ 
is given by
\begin{align} \label{S4-Main-e2}
F_t(z, s)&=z + \sum_{T\in {\mathbb T}^W} t^{|T|}\,  {\Omega} (T, -s)\, \cP_T(z)\\
&=z + \sum_{T\in {\mathbb T}^W}  (-1)^{v(T)} t^{|T|}\,  
\bar{\Omega} (T, s)\, \cP_T(z),\nno 
\end{align}
where $\bar{\Omega} (T, s)$ is the strict 
order polynomial of $T$ as an unlabeled rooted tree. 
\end{theo}

Note that, by the definition of $\theta_T$ 
$(T\in \bbT^W)$ (see page \pageref{theta}), 
we have
\begin{align}
\varphi_T=\theta_{B_+(T)}.\label{varphi-theta}
\end{align}

For some discussions on the order polynomials 
$\Omega (T, s)$ and the strict order polynomials 
$\bar{\Omega} (T, s)$, see \cite{St1} or 
Section $5$ in \cite{GTS-IV}. 
We will need the following results 
on the (strict) order polynomials 
in the proof of the theorem above.

\begin{propo}
For any rooted tree $T$, we have
\begin{align}
\bar \Omega(T, s)&=(-1)^{v(T)} \Omega(T, -s), \label{Recip-Re-e1} \\
\nabla \Omega(T, s)&=\Omega(B_-(T), s), \label{nabla}
\end{align}
where $\nabla: K[s]\to K[s]$ is the linear operator 
that maps any $f(s)\in K[s]$ to $f(s)-f(s-1)$. 
\end{propo}

Eq.\,$(\ref{Recip-Re-e1})$ is a special case of the well-known 
{\it Reciprocity Relation} of the strict order polynomials and order polynomials
of finite posets. For a proof of this remarkable 
result, see Corollary $4.5.15$ in \cite{St1}. 
Eq.\,$(\ref{nabla})$
was first proved by J. Shareshian (unpublished). 
It can be proved by using the definition 
of the strict order polynomials. 
For a proof of a similar property 
of the order polynomials, see Theorem $4.5$ in \cite{WZ}. 
Eq.\,(\ref{nabla}) can be proved by a similar argument as 
the proof there. For more studies on these properties of 
the (strict) order polynomials, 
see \cite{Tree-Inv} and \cite{SWZ}.

\vskip2mm

\underline{Proof of Theorem \ref{S4-Main}}:
 $(a)$  \, By Eq.\,$(\ref{d(t)-DLog})$ and Theorem   
 $(\ref{S3-Main-2})$ 
 and $(\ref{Def-3d(t)})$, 
we have
\begin{align*}
-\lb a_t(z) \pz \rb &=d(t)\\
&={\mathcal A}(\tilde d(t)) \\
&= \sum_{T\in \bar{\mathbb P}^W} t^{|T|} \theta_T {\mathcal D_T}.
\end{align*}

Therefore, by the equation above and Eq.\,$(\ref{D-P})$, we have
\begin{align*}
a_t(z) &=\sum_{T\in \bar{\mathbb P}^W} t^{|T|-1} \theta_T {\mathcal D_T} \cdot z \\
&=\sum_{T\in \bar{\mathbb P}^W} t^{|T|-1} \theta_T \mathcal P_{B_-(T)}(z) \\
\intertext{Replacing the summation index $T\in \bar{\mathbb P}^W$ 
by $B_+(T)$ with $T\in \mathbb T^W$ and noting that $|T|=|B_-(T)|$ 
for any $T\in \bar{\mathbb P}^W$:}
&=\sum_{T\in {\mathbb T}^W} t^{|T|-1} \theta_{B_+(T)} {\mathcal P_T}(z) \\
\intertext{Applying Eq.\,$(\ref{varphi-theta})$:}
&=\sum_{T\in \bT^W} t^{|T|-1} \varphi_T {\mathcal P_T}(z)\,.
\end{align*}
Hence, we get Eq.\,$(\ref{S4-Main-e1})$.

$(b)$ First, let us consider 
the exponential  $e^{s\tilde d(t)}\in 
\cH_{GL}^W [[t]]$. 

\allowdisplaybreaks{
\begin{align}
e^{s \tilde d(t)} &=\sum_{m\geq 0} \frac{s^{m}}{m!}\tilde d(t)^{m} \nno\\
&=1 + \sum_{m\geq 1} \frac{s^{m}}{m!} 
\left( \sum_{T\in \bar{\mathbb P}^W} t^{|T|} \theta_{T} \mathcal V_{T}
   \right)^{m} \nno \\
 \intertext{Applying Lemma \ref{key-lemma-2}:}
   &=1+\sum_{m\geq 1}\frac{s^{m}}{m!} \sum_{T\in \bar{\mathbb T}^W }\,\,\,
  t^{|T|} \sum_{\substack{\vec e=(e_1,\ldots,e_{m-1})\in E(T)^{m-1}\\e_1\succ
\cdots\succ e_{m-1}}}
  \theta_{T_{\vec{e},1}}\cdots
   \theta_{T_{\vec{e}, m}} \mathcal V_T  \nno\\
   &=1 + \sum_{T\in \bar{\mathbb T}^W} \,\left(\sum_{m=1}^{v(T)}\frac{s^{m}}{m!} \,\,\,
   \sum_{\substack{\vec e=(e_1,\ldots,e_{m-1})\in E(T)^{m-1}\\e_1\succ \cdots\succ e_{m-1}}}
\theta_{T_{\vec{e},1}}\cdots
   \theta_{T_{\vec{e}, m}}\right) t^{|T|} \mathcal V_T. \nno
\end{align}} 

Now we apply $\mathcal A$ to the equation above. 
Note that $\mathcal A$ maps $\exp ( s\tilde d(t))$ 
to $\exp(s d(t))$ since, by Theorem \ref{S3-Main-1}, 
$\cA: \cH_{GL}^W \to \cDazz$ is 
a $K$-algebra homomorphism. 
So we have
\begin{align*}
e^{s d(t)}
=1 + \sum_{T\in \bar{\mathbb T}^W} \,\left(\sum_{m=1}^{v(T)}\frac{s^{m}}{m!} \,\,\,
   \sum_{\substack{\vec e=(e_1,\ldots,e_{m-1})\in E(T)^{m-1}\\e_1\succ \cdots\succ e_{m-1}}}
\theta_{T_{\vec{e},1}}\cdots
   \theta_{T_{\vec{e}, m}}\right) t^{|T|} \mathcal D_T. \nno
\end{align*} 

Applying the equation above to $z$ and noting that 
$\mathcal D_T\cdot z= \mathcal P_{B_-(T)}(z)$ 
if $T \in \bar{\mathbb P}^W$ and $0$ otherwise, we get
\begin{align*}
e^{s d(t)} z
&=1 + 
\sum_{T\in \bar{\mathbb P}^W} \,\left(\sum_{m=1}^{v(T)}\frac{s^{m}}{m!} \,\,\,
   \sum_{\substack{e_1\succ \cdots\succ e_{m-1}\\ e_i\in E(T) }}
\theta_{T_{\vec{e},1}}\cdots
   \theta_{T_{\vec{e}, m}}\right) t^{|T|} \mathcal P_T (z) \\
\intertext{Applying Eq.\,$(\ref{theta-orderP-e1})$ in Proposition \ref{theta-orderP}:}
&=1 + \sum_{T\in \bar{\mathbb P}^W}  
\nabla  \Omega (T, s) \mathcal P_T(z).
\end{align*}

Now, replacing $s$ by $-s$ in equation above, 
by Eq.\,$(\ref{Def-flow})$, we get 

\begin{align}
F_t(z, s)&=z + \sum_{T\in \bar{\mathbb P}^W}  \nabla  \Omega (T, -s) \,\,
t^{|T|} \mathcal P_T(z) \nno\\
\intertext{Applying Eq.\,(\ref{nabla}):}
&=z + \sum_{T\in {\mathbb P}^W}  \Omega (B_-(T), -s) t^{|T|} \cP_{B_-(T)} \, (z). \nno\\
\intertext{Changing the summation index $T\in \bar{\mathbb P}^W$ 
by $B_+(T)$ with $T\in \mathbb T^W$:}
&=z + \sum_{T\in {\mathbb T}^W}   \Omega (T, -s) t^{|T|} \cP_T(z) \nno \\
\intertext{Applying Eq.\,$(\ref{Recip-Re-e1})$:}
&=z + \sum_{T\in {\mathbb T}^W}  (-1)^{v(T)} t^{|T|} \Omega (T, s) \cP_T(z).
\end{align}
Hence, we get Eq.\,$(\ref{S4-Main-e2})$.
\epfv

As we mentioned early (see page \pageref{composingPower}),
for any $m\in \bZ$, $F_t(z, m)$ is the $m^{th}$ (composing)
power denoted by $F_t^{[m]}$ of the automorphism $F_t(z)\in \ataz$. 
Hence, by plugging in $m$ for $s$ 
in Eq.\,(\ref{S4-Main-e2}), we get the following formulas.

\begin{corol}\label{M-powers}
For any $\alpha \geq 1$, $m\in \bZ$ and $F_t(z)\in \ataz$, we have
\begin{align} \label{M-powers-e1}
F_t^{[m]}(z)= \sum_{T\in {\mathbb T}^W} t^{|T|} \, {\Omega} (T, -m) \, \cP_T(z),
\end{align}

In particular, by letting $m=-1$, we get the following 
tree expansion formula for the inverse map $G_t(z)$ of $F_t(z)$. 
\begin{align} \label{M-powers-e2}
G_t(z)= \sum_{T\in {\mathbb T}^W}  t^{|T|} \, \cP_T(z),
\end{align}
\end{corol}

\pf Note that, we only need prove Eq.\,(\ref{M-powers-e2}).
But it follows from Eq.\,(\ref{M-powers-e2}) 
and the well-known fact 
that $\Omega(P, 1)=1$ for any finite posets.
\epfv

One remark on the tree 
expansion formulas derived in 
Theorem \ref{S4-Main} 
and Corollary \ref{M-powers} is as follows. 
Note that, one of the conditions 
we have required on $F_t(z)=z-H_t(z)$ 
is that $H_{t=0}(z)=0$. But for 
the automorphisms $F(z)$ of $\kzz$ of the form
$F(z)=z-H(z)$ with $H(z)\in \kzz^{\times n}$ 
and $o(H(z))\geq \alpha$, 
all formulas derived can still be applied to 
$F(z)$ as follows. 

First, we consider the deformation 
$F_t(z)=z-tH(z)$ which does lie in $\ataz$.
Actually, it can be viewed as 
an automorphism of $\ktzz$, 
instead of $\kttzz$, 
over the polynomial algebra $K[t]$.
Therefore, all the formulas above 
with $W=\{1\}$ still apply to $F_t(z)$.
Secondly, by the fact that 
$H_t(z) \in \ktzz^{\times n}$, 
it is easy to check that, 
for any $t_0\in K$, 
the D-Log $a_{t=t_0}(z)$ of $F_{t=t_0}(z)$ 
and the inverse map $G_{t=t_0}(z)$ 
all make sense. In particular, 
by setting $t=1$, we recover 
the D-Log and the formal flow of 
the original automorphism $F(z)$.
Thirdly, in the case $W=\{1\}$, 
the weight $|T|$ $(T\in \bT^W)$ 
of $T$ is same as the number $v(T)$ of 
the vertices of $T$, and the set $\bT^W$ of 
$W$-labeled trees 
can be identify with the set of 
unlabeled rooted trees $\bT$.

By the discussions above, it is easy to see that 
we have the following  tree expansion formulas for 
the automorphisms of $\kzz$.

\begin{corol}\label{S4-Main-C}
For any automorphism $F(z)$ of 
the form $F(z)=z-H(z)$ 
with $o(H(z))\geq 2$, we have

$(a)$ The D-Log $a(z)$ of $F(z)$ is given by 
\begin{align}  \label{S4-Main-C-e1}
a(z)= \sum_{T\in {\mathbb T}} \varphi_T {\mathcal P_T}(z), 
\end{align}
where,  $\varphi_T$ is the 
coefficient of $s$ in the order polynomial 
$\Omega (T, s)$ of the rooted tree $T$.

$(b)$ The formal flow $F(z, s)$ generated by $F(z)$ 
is given by
\begin{align} \label{S4-Main-C-e2}
F(z, s)=z + \sum_{T\in {\mathbb T}}  (-1)^{v(T)} \,  \Omega (T, s)\, \cP_T(z)
\end{align}
where $ \bar{\Omega} (T, s)$ is the strict 
order polynomial of $T$ as an unlabeled rooted tree. 

In particular, we have the following tree expansion inversion formula: 
\begin{align} \label{S4-Main-C-e3}
G (z):=F^{-1}(z)= \sum_{T\in {\mathbb T} } \, \cP_T(z),
\end{align}
\end{corol}

\renewcommand{\theequation}{\thesection.\arabic{equation}}
\renewcommand{\therema}{\thesection.\arabic{rema}}
\setcounter{equation}{0}
\setcounter{rema}{0}

\section[{\bf More Properties of the Specializations 
$\cS$ and $\cT$}]
{\bf More Properties of the Specializations 
$\cS$ and $\cT$ of NCSF's}\label{S5}

Let $K$, $z$, $t$, $\alpha\geq 1$ and $W\subseteq \bN^+$ 
as before. In this subsection, we study more properties of 
the specializations 
$\cS_{F_t}: \cNsf \to \cDazz$ ($F_t\in \ataz$)
in Theorem \ref{S-Correspondence-2} and $\cT_W: \cNsf \to \cH_{GL}^W$ 
in Theorem \ref{Main-Thm-Trees}. We first 
show in Theorem \ref{NCSF-Injc-Trees} that, 
when $W=\bN^+$, the specialization 
$\cT: \cNsf \to \cH_{GL}^W$ of NCSF's 
is actually an embedding. Then, in Theorem \ref{StabInjc-best}, 
we use Theorem $4.6$ in \cite{GTS-II} and
improve it to the family of 
the specializations $\cS_{F_t}: \cNsf \to \cDazz$ 
with all $n\geq 1$ and $F_t= z-H_t(z) \in \ataz$ 
such that $H_t(z)$ is homogeneous and 
the Jacobian matrix $JH$ is strictly lower triangular.  
 
Let us start with the following theorem.

\begin{theo}\label{NCSF-Injc-Trees}
When $W=\bN^+$, the graded $K$-Hopf algebra 
homomorphism $\cT_W: \cNsf \to {\mathcal H}_{GL}^W$ 
in Theorem \ref{Main-Thm-Trees} is an embedding.  
\end{theo}

\pf Let $P \in \cNsf$ be any non-zero NCSF. 
By Theorem $4.6$ in \cite{GTS-II}, 
there exist $F_t\in \ataz$ such that $\cS_{F_t}(P) \neq 0$.
By Proposition \ref{CD-2}, we have 
$\mathcal A_{F_t} (\cT_W (P))= \cS_{F_t}(P)\neq 0$. 
Hence $\cT_W (P)\neq 0$.
\epfv

Combining Corollary \ref{S-cQsf} and
Theorem \ref{NCSF-Injc-Trees} above, 
we get the following corollary.

\begin{corol}\label{Onto}
For any non-empty $W\subseteq\bN^+$, let $\cT:\cNsf\to \cH_{GL}^W$ 
be the specialization of NCSF in Theorem \ref{Main-Thm-Trees} and
$\cT^*$ its graded dual map. 
Then  $\cT^*:\cH_{CK}^W\to {\mathcal Q}sym$ is a homomorphism of 
graded $K$-Hopf algebras.  Furthermore, when $W=\bN^+$, $\cT^*$ is also onto.
\end{corol}

To formulate next main theorem of this section. 
Let us first introduce the following notations.

For any $z$ and $\alpha\geq 1$ as before, 
we let $\mathbb B^{[\alpha]}_t\langle z \rangle$
be the set of automorphisms $F_t=z-H_t(z)$ 
of the polynomial algebra 
$\ktz$ over $K[t]$ such that 
the following conditions are satisfied.
\begin{enumerate}
\item[$\bullet$]  $H_{t=0}(z)=0$. 
\item[$\bullet$] 
$H_t(z)$ is homogeneous in $z$ 
of degree $d \geq \alpha$.
\item[$\bullet$] 
With a proper permutation of the free variables $z_i$'s,
the Jacobian matrix $JH_t(z)$ becomes strictly lower triangular.
\end{enumerate}

Our next main result of this section 
is the following theorem.

\begin{theo}\label{StabInjc-best}
In both commutative and noncommutative cases, 
the following statement holds. 
\vskip2mm
For any fixed $\alpha\geq 1$ and non-zero $P \in \cNsf$, 
there exist $n\geq 1$ 
$($the number  of the free variable 
$z_i$'s$)$ and $F_t(z)\in \mathbb B^{[\alpha]}_t\langle z\rangle$ 
such that $\cS_{F_t} (P)\neq 0$.
\end{theo}

To prove the theorem above, 
we first need the following lemma, 
which is also interesting in its own right.

\begin{lemma}\label{Stab-Injc-trees}
Let $\alpha \geq 1$ and 
$W \subseteq \bN^+$ be fixed above.
For any $T \in \bT^W$, 
there exist $n\geq 1$ and $F_t(z) 
\in \mathbb B^{[\alpha]}_t\langle z\rangle$
such that $P_T \neq 0$ and 
$P_{T'}(z)=0$ for any $T'\in \bT^W$ 
with $|T'|\geq |T|$ but $T'\not \simeq T$.
\end{lemma}

Note that, in the commutative case with $W=\{1\}$, 
the lemma is essentially same 
as Theorem $2.4$ in \cite{WZ}. 
The proof given below is 
also parallel to the proof there.

\pf $(a)$ For any fixed $T\in \bT^W$, we construct 
automorphism $F_t\in \mathbb B^{[\alpha]}_t\langle z\rangle$ 
as follows. Let $n=v(T)$, the number of vertices of $T$, 
and $d \geq \alpha$ be a positive integer that 
is greater or equal to the number of children 
of any vertex of $T$.  
Let $z=(z_1, z_2, \cdots, z_n, z_{n+1})$ 
be free variables.
We first label the edges by $e_2, \ldots, e_n$ 
in an order preserving way, i.e. for any $2\leq i <j \leq n$, 
we have $e_j \succ e_i$. 
We then assign the variable $z_i$ ($2 \le i \le n$) 
to the edge $e_i$ and label the vertices as follows: 
let $v_1=\text{rt}_T$, and for $i=2, \ldots, n$ let $v_i$
be the vertex of $e_i$ which is further away from the root.
Finally, we define $H_t(z) \in \kz^{\times n}$ as follows. 
First, for any $1 \leq i \leq n$ and $m\in W$, 
if $v_i$ is not a leaf of $T$, we set 
$\tilde H_{[m], i}(z) \in \kz$ be the product 
in any fixed order of all the free variables 
assigned to the edges connecting $v_i$ with 
its $m$-labeled children, if there are any, 
and $0$ otherwise. Note that 
$\deg \tilde H_{[m], i}(z) \leq d$.
We set $H_{[m], i}(z)=z_{n+1}^k \tilde H_{[m], i}(z)$ 
for some $k\geq 0$ such that $\deg H_{[m], i}(z)=d$.
Now suppose $v_i$ is a leaf of $T$, we simply set
$H_{[m], i}(z):=z_{n+1}^D$ if $m=\min W$ and $0$ otherwise.
Finally, we set $H_{[m], n+1}(z):=0$ 
for all $m\in W$.
Next, we set $H_{t, i} (z):= \sum_{m\in W} t^m H_{[m], i}(z)$ 
for any $1\leq i\leq n+1$
and $H_t(z)=(\, H_{t,1}(z), H_{t,2}(z), \cdots, H_{t, n+1}(z))$.
Then the wanted automorphism $F_t$ associated with the fixed 
$T\in \bbT^W$ will be $F_t(z):=z-H_t(z)$.

From the construction of $H_t(z)$ above 
and the definition of $P_T(z)$ 
in Eq.\,(\ref{Def-PT}), 
it is easy to check the following facts:
\begin{enumerate}
\item[$\bullet$]  $H_{t=0}(z)=0$ and $H_t(z)$ is homogeneous 
in $z$ of degree $d \geq \alpha$.
\item[$\bullet$] For any $1\leq i\leq n+1$, $H_{t, i}(z)$ only depends 
on the free variables
$z_j$ with $j>i$. Hence, the Jacobian matrix $JH_t(z)$ 
is strictly lower triangular.
\item[$\bullet$] For any $T'\in \bT^W$, $P_{T'}(z)\neq 0$ only if 
 either $T'\simeq T$ or there exists an admissible cut $C\in \mathcal C(T)$ 
such that $T'$ is isomorphic to one of connected components 
obtained by cutting off the edges of $C$ from $T$. 
\item[$\bullet$] For the fixed $T\in \bbT^W$ itself, 
$P_T(z)=c \, z_{n+1}^b$ for some $c, b\in \bN^+$.
\end{enumerate}
Actually, with a little bit more effect, 
one can show that the constants $c$ and $b$ above 
are given by $c=\alpha(T)$ and $b= nd-(n-1)$. 
But we do not need these facts in our proof here.
 
From the discussions above, we only need show
the polynomial map $F_t(z)$ defined above 
is indeed a polynomial automorphism 
of $\ktz$ over $K[t]$, 
i.e. its inverse is also a polynomial map. 
But, from the third observation listed above 
and the tree expansion inversion formula 
Eq.\,(\ref{M-powers-e2}), 
it is easy to see that 
$G_t(z)$ is also a polynomial map.
Note that, this also follows from the 
following well-known result 
in the inversion problem, namely, 
any polynomial map $F_t(z)$ 
with the Jacobian matrix 
$JF_t$ lower triangular 
and invertible is a polynomial 
automorphism.
\epfv

Now we can prove Theorem \ref{StabInjc-best} as follows.

\underline{Proof Theorem \ref{StabInjc-best}}:
Let $P\in \cNsf $ be any non-zero NCSF. 
We choose $W=\bN^+$ and 
let $\cT:\cNsf \to\cH_{GL}^W$ be the 
specialization of NCSF in 
Theorem \ref{Main-Thm-Trees}, 
which we have shown in Theorem 
\ref{NCSF-Injc-Trees} 
is an embedding. Therefore $\cT(P)\neq 0$.
We write $\cT(P)$ as
\begin{align}\label{StabInjc-best-pe1}
\cT(P)=\sum_{T\in \bbT^W} c_T \, T,
\end{align}
with $c_T\in K$ being all but finitely many zero.

Let $k_0\geq 1$ be the least positive integer
such that $c_T=0$ for any $T\in \bbT^W$ 
with $|T|<k_0$. We choose and fix $S\in \bbT^W$ 
such that $c_S\neq 0$ and $|S|=k_0$. 
We fix any $m\in W$ and, for any $T\in \bbT^W$,
denote by $T_m$ the $W$-labeled rooted trees 
in $\bT^W$ obtained by (re)labeling 
the root of $T$ by $m$. Note that, for any 
$T, T'\in \bbT^W$, $T\simeq T'$ in $\bbT^W$ 
iff $T_m\simeq T_m'$ in $\bT^W$.
Now we apply by Lemma \ref{Stab-Injc-trees} 
to $S_m\in \bT^W$ and choose $F_t\in \mathbb B$ 
such that $P_{S_m}(z)\neq 0$ 
but $P_{T_m}(z)=0$ for any $T\in \bbT^W$ with 
$|T|\geq |S|=k_0$ and $T\not \simeq S$. 
Now we apply $\cA_{F_t}$ to $\cT(P)$. 
Note that, by Eq.\,(\ref{LL4.1.2-e1}), 
we have, $D_T \cdot H_{[m]}(z)=P_{T_m}$ 
for any $T\in \bbT^W$. 
By using all the facts above and 
Eq.\,(\ref{StabInjc-best-pe1}), 
it is easy to see that, we have
\begin{align*}
\cA_{F_t}(\cT (P)) \cdot H_{[m]}(z)&= 
\sum_{T\in \bbT^W} c_T \, P_{T_m}(z)\\
&= P_{S_m}(z) \neq 0.
\end{align*}

Finally, by Theorem \ref{CD-2}, we have 
$\cS_{F_t}(P)=\cA_{F_t} ( \cT (P))\neq 0$.
Hence, we have proved 
Theorem \ref{StabInjc-best}.
\epfv

{\small \sc Department of Mathematics, Illinois State University,
Normal, IL 61790-4520.}

{\em E-mail}: wzhao@ilstu.edu.

\end{document}